\documentclass{amsart}
\usepackage[utf8]{inputenc}
\usepackage{amsfonts}
\usepackage{hyperref}
\usepackage{amsmath}
\usepackage{xcolor}
\usepackage{color}
\usepackage{amsthm}
\usepackage{pdflscape}
\usepackage{mathrsfs}
\usepackage{verbatim}
\usepackage{mathtools}
\usepackage{dsfont}
\usepackage{bbm}
\usepackage{amssymb}
\usepackage{geometry}
\usepackage{amstext}
\usepackage{enumerate}
\usepackage{graphicx}
\usepackage{subfig}
\usepackage{geometry}

\usepackage{float} % I am trying to fix the images

\usepackage{xcolor}
\usepackage{hyperref}
\hypersetup{colorlinks,% 
citecolor=red,% 
linkcolor=blue,% 
}

\theoremstyle{plain}
\newtheorem*{thm*}{Theorem}
\newtheorem{thm}{Theorem}[]
\newtheorem{cor}{Corollary}[]

\newtheorem{lem}{Lemma}[]
\newtheorem{conj}{Conjecture}[]
\newtheorem{question}{Question}[]

\newtheorem{rem}{Remark}[]
\newtheorem{example}{Example}[]

\newtheorem{observation}{Observation}[]

\theoremstyle{definition}
\newtheorem{defn}{Definition}

\newcommand{\ord}{{\rm ord}}
\newcommand{\tree}{{\rm Tree}}
\newcommand{\N}{\mathbb{N}}
\newcommand{\Z}{\mathbb{Z}}

\newcommand{\C}{\mathbb{C}}

\newcommand{\RR}{\mathcal{R}}

\newcommand{\blue}[1]{\textcolor{blue}{#1}}

\newcommand{\Mod}[1]{\ \left(\mathrm{mod}\ #1\right)}

\title{Rational eigenfunctions of the Hecke operators}

\author{Andr\'e Rosenbaum Coelho}
\address{Instituto de Matem\'atica e Estat\'\i stica, Universidade de S\~ao Paulo\\ Rua do Mat\~ao 1010, 05508-090 S\~ao Paulo/SP, Brazil.}
\email{andrerosenbaumcoelho@usp.br}

\author{Caio Simon de Oliveira}
\thanks{The second author is grateful for the support of CNPq 2023-2684}
\address{Instituto de Matem\'atica e Estat\'\i stica, Universidade de S\~ao Paulo\\ Rua do Mat\~ao 1010, 05508-090 S\~ao Paulo/SP, Brazil.}
\email{caiosimon@usp.br}

\author{Sinai Robins}
\address{Instituto de Matem\'atica e Estat\'\i stica, Universidade de S\~ao Paulo\\ Rua do Mat\~ao 1010, 05508-090 S\~ao Paulo/SP, Brazil.}
\email{srobins@ime.usp.br}

\begin{document}

\begin{abstract}
We study the action of the Hecke operators $U_n$ on the space $\RR$ of rational functions in one variable, over $\C$. 
The main goal is to give a complete classification of the  eigenfunctions of $U_n$.
We accomplish this by introducing certain number-theoretic directed graphs, called Zolotarev Graphs, which extend the well-known permutations due to Zolotarev.

We develop the theory of these Zolotarev graphs, using them to decompose the eigenfunctions of $U_n$ into certain natural finite-dimensional vector spaces of rational functions, which we call the eigenspaces.  In this context, we prove that the dimension of each eigenspace is equal to the number of nodes of a cycle that belongs to its corresponding Zolotarev graph. We prove that the number of leaves of this Zolotarev graph equals the dimension of the kernel of $U_n$.  We then give a novel number-theoretic formula for the number of cycles of fixed length, in each Zolotarev graph.  We also study the simultaneous eigenfunctions  for all of the $U_n$, and give explicit bases for all of them.  In the process, we answer many questions that were set out in \cite{GilRobins}.  We also discover certain strong relations between these graphs and the kernel of $U_n$ acting on a subspace of $\RR$; in particular, we give several equivalent conditions for the diagonalizibility of $U_n$.

Finally, we prove that the classical Artin Conjecture on primitive roots is equivalent to a new conjecture here, that infinitely many of these eigenspaces have dimension $1$.  
\end{abstract}
\maketitle

\tableofcontents

%%%%%%%%%%%%%%%%%%%%%%%%%%%
%%%%%%%%%%%%%%%%%%%%%%%%%%%

\section{Introduction} \label{sec:intro}
A very useful linear operator in Number theory is the series dissection operator:
\begin{equation} \label{def:Hecke operator}
U_n\left(\sum_{k=0}^{\infty} a(k) x^k\right):= \sum_{k=0}^{\infty} a(nk) x^k,
\end{equation}
where $n$ is any fixed positive integer. 
The operator $U_n$ was thoroughly studied  by Hecke in the important case that
 the series in \eqref{def:Hecke operator} represents a modular form, and is therefore called a {\bf Hecke operator}  (see \cite{Hecke} for Erich Hecke´s original paper).  In other contexts such operators are also called dissection operators (they pick out every $n$'th coefficient of the series). 
 In Vlasenko's work \cite{Beukers.Vlasenko}, the Hecke operator $U_p$ is also called a Cartier operator.
 
 Given a positive integer $n$, here we study the action of $U_n$ on the vector space of rational functions over $\C$ that are analytic at the origin:
 \begin{equation}
 \RR:= \left  \{ 
 \frac{p(x)}{q(x)} \mid 
 p, q \in \C[x] \text{ and } \frac{p(x)}{q(x)}= \sum_{k=0}^\infty a(k) x^k,
 \text{ near the origin}
 \right \}.
 \end{equation}

The operator
$U_n$ acts on the vector space $\RR$, sending a rational function to another rational function (see \cite{GilRobins}).  It follows easily from its definition that $U_n$ is a linear operator.

Our first main goal is to give an explicit description of all the rational functions in $\RR$ that are eigenfunctions of $U_n$:
\begin{equation}
\label{def:eigenfunction of U_n}
    U_n f(x) = \lambda f(x),
\end{equation}
for all $x$ near the origin.
For each $n \in \Z_{>0}$, there are always rational eigenfunctions of $U_n$, for example the trivial eigenfunction $f(x):= \frac{1}{1-x}= \sum_{k\geq 0} x^n$, which satisfies \eqref{def:eigenfunction of U_n}.
In \cite{GilRobins}, we studied and classified the eigenvalues of the Hecke operators $U_n$, on the space of rational functions with real coefficients, but \cite{GilRobins} did not classify the eigenfunctions.    Here we extend this study by completely classifying all of the eigenfunctions of $U_n$, as well as working more generally over $\C$.  However, the methods of this paper are essentially independent of the methods of \cite{GilRobins}, and indeed we proceed from first principles.

Our second main goal is to introduce and study the related number-theoretic graphs, called Zolotarev graphs (Section \ref{sec: Zolotarev graphs}), which we use to classify the eigenfunctions of $U_n$.

Our final objective is to restate the Artin Conjecture on primitive roots in terms of one-dimensional eigenspaces of $U_n$.

%A crucial idea here (was) to iterate the linear operator $U_n$ a finite number of times in order to  study the orbits of certain rational functions, an instance of  \emph{arithmetic dynamics} (\cite{Silverman.book}, \cite{Benedetto.etal}). \orange{Probably not so crucial anymore...}\orange{Yes...}
 
%%%%%%%%%%%%%%%%%%%%%%%%%%%%%%%%%%%%
\medskip
\begin{thm}
\label{thm:first result}
Fix a positive integer $n>1$, and 
let $f\in \RR$ be an eigenfunction of $U_n$ with eigenvalue $\lambda\neq0$. 
Then all of the poles of $f$ are roots of unity.

\hfill
\hyperlink{proof of thm:first result}{(Proof)}
 $\square$
\end{thm}
By  Theorem \ref{thm:first result}, there is a  
{\bf least positive integer $L$} such that the poles of $f$ are all $L$'th roots of unity. We call  this unique $L$ the {\bf level} of $f$, and we will sometimes use the notation $\rm{level}(f) = L$. 

\begin{lem}
\label{gcd(n, L)=1 for an eigenfunction}
Suppose $f\in \RR$ is an eigenfunction of $U_n$, with $n>1$, such that $f$ has level $L$, and eigenvalue $\lambda\neq0$. 
Then $\gcd(n, L)=1$.

\hfill 
\hyperlink{proof of gcd(n, L)=1 for an eigenfunction}{(Proof)}
$\square$
\end{lem}

%\begin{lem}
%\label{lemma:How to write the Taylor Series of an eigenfunction}
%Suppose $f\in \RR$ is an eigenfunction of $U_n$ with $n>1$ and eigenvalue $\lambda\neq0$. Then the Taylor series coefficients $a(k)$ of $f$ have the form:
%\[
%a(k)=\sum_{j=1}^{d}C_jk^{m_j-1}r_j^{k}
%\]
%where $m_j$ is the multiplicity of each distinct pole $\frac{1}{r_j}$.
%
%\hfill 
%\hyperlink{proof of thm:second result}{(Proof)}
%$\square$
%\end{lem}

\begin{thm}[The weight of an eigenfunction]
\label{thm:second result}
Fix a positive integer $n>1$, and 
let $f\in \RR$ be an eigenfunction of $U_n$, with eigenvalue $\lambda \not=0$. Then all
 the poles of $f$ have the same multiplicity, which we denote by $\kappa$.
 
 \hfill
 \hyperlink{proof of thm:second result}{(Proof)}
 $\square$
\end{thm}
%%%%%%%%%%%%%%%%%%
\noindent
 Moreover, under the assumption of 
Theorem \ref{thm:second result}, 
we call the unique nonnegative integer $\kappa$ 
the {\bf weight} of $f$, and write 
$\rm{weight}(f)=\kappa$.
We first observe an interesting relation between $n$ and the level $L$ of an eigenfunction.

\noindent
Now it makes sense to work with the order of $n$ modulo $L$, which we denote throughout by 
\[
c:= {\rm ord}_L(n).
\]
\begin{thm}
\label{thm:eigenvalues and kappa}
Fix a positive integer $n>1$, and 
let $f\in \RR$ be an eigenfunction of $U_n$.
If we have $U_n f(x) = \lambda f(x)$, for a nonzero eigenvalue $\lambda$, then 
\[
\lambda = \omega \, n^{\kappa-1}, 
\]
where 
\[
\omega:= e^{ \frac{2\pi i m}{c}},
\]
for some integer $m$.  
Here $\rm{weight}(f)=\kappa$, and
${\rm level}(f)=L$.

\hfill
 \hyperlink{proof of thm:second result}{(Proof)}
 $\square$
\end{thm}

%%%%%%%%%%%%%%%%%%%%%%%%%%%%%%%%%%%%
\begin{cor}
\label{corollary:How to write an eigenfunction}
For each integer $n>1$, the eigenfunctions $f(x)$ of $U_n$ with eigenvalue $\lambda\neq0$, level $L$ and weight $\kappa$ may be written as
\[
    f(x)=\sum_{k=0}^{\infty} k^{\kappa-1}a(k)x^k
\]
where $a$ is a periodic function on $\Z_{\geq 0}$, with smallest period $L$.
\hfill 
\hyperlink{proof of thm:second result}{(Proof)}
$\square$
\end{cor}
%%%%%%%%%%%%%%%%%%%%%%%%%%%%%%%%%%%%
\medskip
\begin{cor}[The spectrum of $U_n$]
\label{cor:the spectrum of U_n}
    The spectrum of the operator $U_n$, acting on $\RR$, is:
\begin{equation}
{\rm Spec}(U_n)=  
\{    
e^{ \frac{2\pi i m}{c}} n^{\kappa-1}  \mid
m \in \Z, \text{ and } \kappa, L \text{ are  positive integers}
, (n,L)=1
\},
\end{equation}
where $c:=\ord_L(n)$.
\hfill 
\hyperlink{proof of cor:the spectrum of U_n}{(Proof)}
$\square$
\end{cor}

\begin{cor}
\label{cor: image of phi function}
    If $N$ lies in the image of Euler's totient function $\phi$, then $\lambda=e^{\frac{2\pi i}{N}}$ is an eigenvalue of $U_n$.
\hfill
\hyperlink{proof of cor: image of phi function}{(Proof)}
$\square$
\end{cor}
\noindent
The converse of Corollary \ref{cor: image of phi function}
is false (see Remark \eqref{converse of cor. 3 is false}).

%%%%%%%%%%%%%%%%%%%%%%%%%%%%%%%%%%%%%%%%%%%%%%%%%%%%%%

\noindent 
Corollary \ref{corollary:How to write an eigenfunction} strongly suggests that we consider the following natural subspace of $\RR$: 
\begin{equation}
  \RR(L, \kappa):= \left\{
  f\in \RR \mid 
  f(x)=\sum_{k\geq0}k^{\kappa-1}a(k)x^k, 
  \text{ where } a(k+L) = a(k)
  \right\},
\end{equation}
which also contains functions that are not eigenfunctions of any $U_n$. We note that if $L_1\mid L_2$, then $\RR(L_1,\kappa)\subseteq \RR(L_2,\kappa).$
Throughout, we'll use the following useful notation.
\begin{defn}
The indicator function of an arithmetic progression $\{ a + nL \mid n \in \Z\}$ is defined by:
\begin{equation}\label{Indicator Function}
1_{a\text{ mod }L}(x):=
\begin{cases}
0 & \text{ if } x \not\equiv a \Mod L \\
1 & \text{ if } x  \equiv a  \Mod L.
\end{cases}
\end{equation}
\end{defn}

\begin{lem}
\label{lem: precise dimension of R_L,k}
    $\RR(L,\kappa)$ is a vector space with $\dim \RR(L,\kappa)=L.$
    \hfill
    \hyperlink{proof of lem: precise dimension of R_L,k}{(Proof)}
$\square$
\end{lem}

%%%%%%%%%%%%%%%%%%
%%%%%%%%%%%%%%%%%%

\noindent 
We fix two positive integers $\kappa$ and $L$.  For each given root of unity $\omega \in \C$, we consider the eigenspace corresponding to $\omega$: 
\begin{equation}
\label{def: eigenspace}
  E_n(\omega, L, \kappa):= \left\{ f\in \RR \mid 
  U_n(f)  = n^{\kappa-1}\omega f,  
  \text{ where } \text{level}(f) \mid 
L, \text{and weight}(f)= \kappa\right\},
\end{equation}
one of the most important constructions here.
We observe that if $L_1\mid L_2$, then $E_n(\omega,L_1,\kappa)\subseteq E_n(\omega,L_2,\kappa)$, because 
\[
E_n(\omega, L, \kappa)= \left\{ f\in \RR(L,\kappa) \mid 
  U_n(f)  = n^{\kappa-1}\omega f\right\}.
  \]

%where $\left< v\in S\right>$ is by definition the subspace of $\RR$ consisting of all complex linear combinations of vectors belonging to $S\subset \RR$. 

%%%%%%%%%%%%%%%%%%
Each vector space $E_n(\omega, L, \kappa)$ is a convenient packet of eigenfunctions of $U_n$.  It'll turn
out that $E_n(\omega, L, \kappa)$ forms a very interesting
  finite dimensional vector space, whose elements have number-theoretic content.  To describe the rational eigenfunctions of $U_n$, we'll develop some special permutations first.
%%%%%%%%%%%%%%%%%%

\subsection{Zolotarev permutations}
There are some natural permutations that we study here, which naturally arise in the 
study of the eigenfunctions of Hecke operators. 
 Throughout this section, we let $n,L\in\Z_{>0}$,  $\left(n,L\right)=1$, and 
$c:={\rm ord}_L(n)$.
Consider the map that sends each element $m\in \Z/L\Z$ to $mn \pmod L$.
We will write $\overline{mn}$ for the unique integer in the interval $[0, L-1]$ that is congruent to $mn \pmod L$.  Because $(n, L)=1$, this map is a permutation of 
$\Z/L\Z$, which we may write as:
\begin{equation}
\label{modular permutation}
\tau(n, L):= 
\begin{pmatrix}
0 & 1 & \cdots & k & \cdots & L-1\\
 &   &   &   &     & \\
\overline{0} & \overline{n} & \cdots & \overline{kn} & \cdots & \overline{(L-1)n} 
\end{pmatrix},
\end{equation}
which we'll call a {\bf Zolotarev permutation}, because these permutations were used by Zolotarev to prove the reciprocity laws of Gauss and of Jacobi.  

The permutation $\tau(n, L)$ breaks up into a product of disjoint cycles (called the disjoint cycle decomposition) where for each $1\leq j \leq c$, we define \begin{equation}
    b_j:= \text{ the number of disjoint cycles of length } j, \text{ in } \tau(n, L).
\end{equation} 
Of course, we allow $b_j=0$ to hold as well.  We note that there always exists the following cycle in the permutation $\tau(n, L)$:
 $(1 \rightarrow \overline{n} \rightarrow \overline{n}^2 \rightarrow \cdots \overline{n}^{c-1})$.
Zolotarev permutations have also been studied in 
\cite{BennettMosteig}, for example, in the context of cyclotomic cosets.
%\begin{equation}
%\label{def:cyclotomic coset}
%    C_{n, L, m} = \left\{
 % n^i m\mod L \mid  i \geq 1
%\right\}, 
%\text{ where } \gcd(n, L)=1.
%\end{equation}

\begin{example}
\rm{Let $L=30$, and $n=7$.  It's easy to check that 
${\rm ord}_{30}(7)= 4$.  Here \eqref{modular permutation} has the disjoint cycle decomposition:
\[
\tau(7, 30)= (1 \ 7 \ 19 \ 13)
(2 \ 14  \ 8 \ 26)(3 \ 21 \ 27 \ 9)(4 \ 28 \ 16 \ 22)(6 \ 12 \ 24 \ 18)(11 \ 17 \ 29 \ 23)
(0)(5)(10)(15)(20)(25).
\]
Therefore $b_1= 6$, $b_2=b_3=0$, and $b_4 = 6$.
\hfill $\square$}
\end{example}
\noindent
It is natural to ask: does there exist some nice formula for $b_j$? 
\begin{lem}[Formula for the number of cycles]
\label{formula for b_j}
We fix positive integers $n, L$. Then the number of cycles of length $j$ in 
$\tau(n, L)$ is given by
\begin{equation}
\label{Formula for b_j}
b_j = \frac{1}{j}
\sum_{d \mid j} 
\mu\left(\frac{j}{d}\right)
\gcd(n^d-1, L).
\end{equation}
\hfill 
\hyperlink{proof of Lemma: formula for the b_j}{(Proof)}
$\square$
\end{lem}
We observe that the proof of Lemma \ref{formula for b_j} does not require $n$ to be coprime to $L$.   This suggests a more general theory than Zolotarev permutations, which we call Zolotarev graphs, and which form an integral part of our study of eigenfunctions of $U_n$, beginning with Section \ref{sec: Zolotarev graphs}.

\begin{example}
\rm{We've already seen that 
$b_1 = \gcd(n, L-1)$, by Lemma \ref{formula for b_j}.  If $L := p$, a prime, then  all cycles in $\tau(n, L)$ must have a length that divides $\phi(p)= p-1$. 
It follows that when we calculate $b_j$, we must have $j \mid p-1$.
By Lemma \ref{formula for b_j}, we have:
\begin{align}
\label{b_j formula for L= prime}
b_j = \frac{1}{j}
\sum_{d \mid j} 
\mu\left(\frac{j}{d}\right)
\gcd(n^d-1, p),
\end{align}
for each $j \mid p-1$.
To see a special case, we let $L=17, n=2$.  Then 
$j | 16\implies j\in \{1, 2, 4, 8, 16\}$.  Using \eqref{b_j formula for L= prime}, we may compute the nonzero $b_j$: \ 
$b_1 = \gcd(n-1, 17)=1, \ 
b_2 = \tfrac{1}{2}
\sum_{d \mid 2} 
\mu\left(\frac{2}{d}\right)
\gcd(2^d-1, 17) =0, \ 
b_4 =0, b_8  = 2$, and $b_{16}=0$.
\hfill $\square$}
\end{example}

\begin{lem}
\label{lem:defining b_j}
Let  $L, n$ be positive coprime integers, and let  
$c:={\rm ord}_L(n)$.  Then the Zolotarev permutations enjoy the following properties:
\begin{enumerate}[(a)]
\item if $\overline{c}$ is the size of a cycle, then $\overline{c}\mid c$.
\item Let $A$ and $B$ be any two cycles in $\tau(n, L)$. 
Then $A$ and $B$ are distinct $\iff$ $A\cap B=\varnothing$.
\item $\sum_{j=1}^{c}jb_j=L$.
    \end{enumerate}
\hfill
\hyperlink{proof of lem:defining b_j}{(Proof)}
$\square$
\end{lem}
\noindent

\subsection{Remarks}
\begin{enumerate}[(a)]
\item 
We note that the proof of Lemma \ref{formula for b_j} uses the elementary fact that 
$
L\mid k(n^j-1) \iff \frac{L}{\gcd(n^j-1, L)}\mid k.
$
But the latter statement is also equivalent to
$\frac{L}{(L,k)}\mid n^j-1$, which gives us another way to calculate $b_j$, as follows. 
\begin{equation}
\label{another way to compute b_j}
b_j=\frac{1}{j}
\left| \left\{ 1\leq m \leq L \mid 
j = \ord_{\frac{L}{(L,m)}}(n)
\right\}
\right|.
\end{equation}
\item 
Equation \eqref{another way to compute b_j} also gives the following identity:
\[
\sum_{d\mid j}
\mu\left(\frac{j}{d}\right)
\gcd(n^d-1,L)
=\left| \left\{ 1\leq m \leq L \mid 
j = \ord_{\frac{L}{(L,m)}}(n)
\right\}
\right|.
\]
\item 
Suppose we fix an eigenfunction 
 $f \in E_n(\omega, L,\kappa)$, with an
eigenvalue $\lambda= n^{\kappa -1} \omega$, where $\omega$ is a root of unity.  Then using the notation of Theorem \ref{thm:eigenvalues and kappa} we know that:
 \[
 \frac{\lambda}{n^{\kappa -1}} =  \omega = e^{ \frac{2\pi i m/\gcd(m, c)}{c/\gcd(m, c)}}:=  e^{ \frac{2\pi i m/\gcd(m, c)}{h}},
 \]
where $h:= \frac{c}{\gcd(m, c)}=\frac{\ord_L(n)}{\gcd(m, c)}=\ord_L(n^m)$.
\bigskip
\item \label{converse of cor. 3 is false}
The converse to Corollary \ref{cor: image of phi function} does not hold. Indeed, consider $U_2:\RR(29,1)\rightarrow\RR(29,1)$. By means of an easy computation, we verify that $2$ is a primitive root modulo 29, thus $\ord_{29}(2^2)=\ord_{29}(4)=14.$ Therefore, by Corollary \ref{cor:the spectrum of U_n}, $e^{\frac{2\pi i}{14}}\in\text{Spec}(U_2)$.
On the other hand, \cite{Coleman} shows that the equation $\phi(x)=14$
does not admit a solution.
\end{enumerate}
%%%%%%%%%%%%%%%%%%%%%%%%%%%%%%%%%%%%
%%%%%%%%%%%%%%%%%%
\bigskip
\section{Zolotarev graphs}
\label{sec: Zolotarev graphs}

  To completely characterize the eigenspaces given by $U_n$ acting on $\RR(L,\kappa)$, including its kernel, we must extend the Zolotarev permutations to any two integers $n$, $L$, not necessarily coprime. The fact that we do not require $n$ and $L$ to be coprime in the proof of Lemma \ref{formula for b_j} already foreshadows a generalization of the Zolotarev permutations.

We now introduce the \textbf{Zolotarev graphs}. In the followings sections, we will use this strong tool to paint a more complete picture of the eigenspaces of $U_n$, and of its 
kernel.

\begin{defn} 
\label{def: zolotarev graphs}
Let $Z(n, L)$ be the directed graph whose nodes are given by the integers in $\Z/L\Z$, and such that there exists a directed edge from node $a$ to  node $b$  precisely when $n a \equiv b \Mod L$.  We'll often use the notation $a\longmapsto b$ and call such a directed graph a {\bf Zolotarev graph}.
\end{defn}
So the Zolotarev permutations, defined by \eqref{modular permutation}, are now replaced by Zolotarev  graphs above.
Given a Zolotarev graph $Z(n, L)$, we define the following terms, for any $m \in \Z/L\Z$:
\begin{enumerate}
    \item $m$ is called a {\bf leaf} integer if there does \emph{not} exist an integer $k$ such that $nk \equiv m \Mod L$, i.e., the \textbf{indegree} of $m$ is 0.
    \item $m$ is called a {\bf root} integer if there exists an integer $k$ such that $n^k m\equiv m\Mod L$.
    \item $m$ is called a {\bf branch} integer if it is neither a leaf, nor a root.
\end{enumerate}
For simplicity, each integer in $m \in\Z/L\Z$ will simply be called a leaf, a root, or a branch.  From these definitions, it follows that the sum of the number of leaves, roots and branches in $Z(n,L)$ is precisely equal to $L$.
\begin{defn} 
\label{def: directed path}
A {\bf directed path} of length $N$ is a sequence of nodes $m_1, m_2, \cdots, m_N \in Z(n, L)$, such that there is a directed edge between $m_i$ and $m_{i+1}$, for all $i=1, 2, \dots, N-1$.
\end{defn}
We may also consider the {\bf undirected graph $G(n, L)$} that is defined by using the same data of $Z(n,L)$, but simply ignoring the directions on each edge. 
 We consider a {\bf maximal set of nodes} in $G(n, L)$ such that for each $2$ of its nodes $a, b$, we have an undirected path between $a$ and $b$.  The subgraph of $G(n, L)$, induced by these nodes, defines  a connected component of $G(n, L)$.
\begin{defn} 
\label{def: connected component}
 We define a {\bf connected component} of $Z(n, L)$ by taking a connected component of $G(n, L)$, and including the original direction on each of its edges.  
\end{defn}
\begin{example}
\rm{In Figure \ref{newL=20,n=10,14}, we see the Zolotarev graph $Z(10, 20)$, which has just $1$ connected component, $1$ root, $1$ branch, and $18$ leaves.
%We also see the Zolotarev graph $Z(14, 20)$, which has $3$ connected components, $5$ roots, $5$ branches, and $10$ leaves. 
We note that the indegree of a node that is not a leaf in %$Z(14, 20)$ is equal to $\gcd(14, 20) = 2$, and the indegree  of a node that is not a leaf in 
$Z(10, 20)$ is equal to $\gcd(10, 20) = 10$.}
\begin{figure}[H]
\includegraphics[width=4cm]{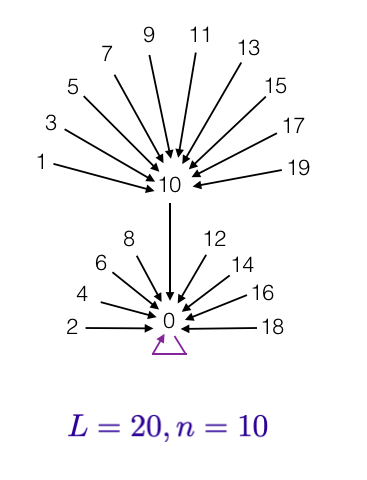}
\caption{The Zolotarev graph $Z(10, 20)$ has $1$ root, $1$ branch, and $18$ leaves.  It has only one connected component. } 
\label{newL=20,n=10,14}
\end{figure}

\end{example}
\noindent
See Appendix \ref{appendix: examples} for all of the Zolotarev Graphs $Z(n, 20)$, with $\gcd(n, 20) >1$. 
Throughout, we will also use the following important function of $L$ and $n$.
\begin{defn}[The simplified $L$, with respect to $n$]
\label{def:L_n}
For any two positive integers $L, n$, we let
\begin{equation}
L_n := \text{ the largest divisor of } L 
\text{ that is coprime with } n.
\end{equation}
\end{defn}
%%%%%%%%%%%%%%%%%%%%%%%%%%%%%%%%%%%%%%

\begin{example}
\label{recursive formula to L_n}
 \rm{Let $L=40$, and $n=4$. 
 Here $L_n= 5$.  We note that there is an alternative algorithm to compute $L_n$. Namely, we may define $x_1:=\frac{L}{\gcd(L, n)}$, and inductively define 
 $x_j:=\frac{x_{j-1}}{\gcd(x_{j-1}, n)}$. Then it's an easy exercise to show that after a finite number of steps this sequence converges to $L_n$.  In this example, we have $x_1 = 10$, and $x_2 = 5= L_n$.}
\hfill $\square$
\end{example}

\medskip

%%%%%%%%%%%%%%%%%%%%%%%%%%%%%%%%%%

\subsection{Structure of the Zolotarev graphs}
\begin{defn}
    Given a Zolotarev Graph $Z(n,L)$ and a root $m$, a \textbf{cycle} of length $j$, containing $m$ (if it exists), is the set 
\[
    A=\{  m\cdot n^k \mid  k\in\Z  \}, \text{ with } |A|=j.
\]
\end{defn}

We note that it follows easily from the definitions above that each connected component of $Z(n, L)$ contains exactly one cycle. 
\begin{observation}
\label{Zolotarev graphs:formula for b_j}
Let $n, L$ be positive integers. Then the number of cycles of length $j$ in 
$Z(n, L)$ is given by
\begin{equation}
\label{formula for b_j, Zolotarev graphs}
b_j = \frac{1}{j}
\sum_{d \mid j} 
\mu\left(\frac{j}{d}\right)
\gcd(n^d-1, L).
\end{equation}
\end{observation}
We remark that the proof of this observation follows directly from the proof of Lemma \ref{formula for b_j}, since it does not rely on any assumption regarding the greatest common divisor of $n$ and $L$.
\begin{observation}[Number of connected components]
    There are
    \[
   \sum_{j=1}^{L}
\frac{1}{j}\sum_{d \mid j} 
\mu\left(\frac{j}{d}\right)
\gcd(n^d-1, L).  
\]
connected components in $Z(n,L)$, which follows directly from 
\eqref{formula for b_j, Zolotarev graphs}.
\hfill $\square$
\end{observation}

We also observe that all nodes that are not leaves have the same indegree, namely $(n, L)$.

\begin{example}
\rm{Let's compute $b_j$, the number of cycles of length $j$, in $Z(2, 20)$.  By Lemma 
\ref{Zolotarev graphs:formula for b_j}, we compute the following $b_j$ values:
$b_1 = \gcd(2-1, 20)=1$, 
$b_2 = 
\frac{1}{2}
\sum_{d \mid 2} 
\mu\left(\frac{2}{d}\right)
\gcd(2^d - 1, 20)=0$, 
\[
b_4 = 
\frac{1}{4}
\sum_{d \mid 4} 
\mu\left(\frac{4}{d}\right)
\gcd(2^d - 1, 20)= 
\frac{1}{4}
\left( 
-\gcd(2^2 - 1, 20)
+
\gcd(2^4 - 1, 20)
\right)
=1.
\]
We also have $b_j=0$ for all other values of $j$, and the latter data confirms our data from Figure \ref{newL=20,n=2,18}.}
\hfill $\square$
\end{example}

\noindent
We will make repeated use of the following  useful fact.
\begin{lem}
\label{lem: (n,L) arrows}
Suppose that $k\in Z(n,L)$ is not a leaf, 
and let $g:= \gcd(n, L)$.   
     Then the $g$ solutions to $n\cdot c \equiv k \Mod L$ are all congruent to a fixed $c$ mod $\frac{L}{g}$.
\hfill
\hyperlink{proof of lem: (n,L) arrows}{(Proof)}
$\square$
\end{lem}
In fact, our Algorithm Z below 
(section \ref{Algorithm:Zolotarev graphs}) makes repeated use of Lemma \ref{lem: (n,L) arrows}, by commencing with a smaller modulus, and building up the corresponding  inequivalent solutions for the higher modulus.

%\begin{observation}
%\label{obs: connectivity of Z(n,L)}
 %  Every node in $Z(n,L)$ is connected to a root.
%\end{observation}

\begin{thm}[Arithmetic equivalences for leaves, branches, and roots]
\label{thm: zolotarev structure}
If $(L,n)>1$ and $Z(n,L)$ is a Zolotarev graph, the following hold:
    \begin{enumerate}
        \item $m$ is a leaf $\iff$ $(n,L)\nmid m$
        \item $m$ is a root $\iff$ $\frac{L}{L_n}\mid m$
        \item $m$ is a branch $\iff$ $(n,L)\mid m$ and $\frac{L}{L_n}\nmid m$
    \end{enumerate}
\hfill
\hyperlink{proof of thm: zolotarev structure}{(Proof)}
 $\square$
\end{thm}

\begin{cor}[Number of roots, leaves, and branches]
\label{cor: number of roots, leaves and branches}
Any Zolotarev graph $Z(n, L)$ has precisely:
\begin{enumerate}
\item $L_n$ roots.
\item $L-\frac{L}{(n,L)}=L\frac{(n,L)-1}{(n,L)}$ leaves.
\item $\frac{L}{(n,L)}-L_n$ branches.
\end{enumerate}
\hfill
\hyperlink{proof of cor: number of roots, leaves and branches}{(Proof)}
$\square$
\end{cor}
We'll show in Section \ref{sec: The Kernel}
below that the  dimension of the kernel of $U_n$ acting on $\RR(L,\kappa)$ is equal to the number of leaves of $Z(n,L).$   

\begin{lem}
\label{lem: m+kL_n is a root?}
    If $m\in Z(n,L)$, then exists a unique integer $k$, with $0\leq k <\frac{L}{L_n}$ such that $m+kL_n$ is a root.
\hfill
\hyperlink{proof of lem: m+kL_n is a root?}{(Proof)}
$\square$
\end{lem}

\begin{thm}[Alternate formula for the number of roots]
\label{thm: number of roots with b_j}
The number of roots of a Zolotarev graph 
$Z(n, L)$ is equal to $\sum_{j=1}^{L}{jb_j}$, and may be expressed as follows:
\begin{equation}
|roots|= 
 \sum_{j=1}^{L}
\sum_{d \mid j} 
\mu\left(\frac{j}{d}\right)
\gcd(n^d-1, L).
\end{equation}
\hfill
\hyperlink{proof of thm: number of roots with b_j}{(Proof)}
$\square$
\end{thm}

\noindent
We immediately obtain from Corollary \ref{cor: number of roots, leaves and branches} and Theorem \ref{thm: number of roots with b_j} an explicit formula  for $L_n$, as follows.
\begin{cor}
We have the following explicit formula for $L_n$:
\begin{equation}
\label{eq: explicit formula to L_n}
    L_n  = 
   \sum_{j=1}^{L}
\sum_{d \mid j} 
\mu\left(\frac{j}{d}\right)
\gcd(n^d-1, L).  
\end{equation}
\hfill $\square$ 
\end{cor}

\begin{defn}
  We define the {\bf distance between any two nodes} $a,b$ of $Z(n, L)$, 
to be the smallest integer $k$ that satisfies $n^k a\equiv b\mod L$, if such a $k$ exists and we write $d(a,b)=k$.
In case the latter congruence does not have a solution, we say that $d(a,b) = \infty$.
\end{defn}
In particular, we note that $d(a,a) = 0$, for all $a \in Z(n, L)$.

\begin{defn}
    Given any root $r$ of $Z(n, L)$, we consider the directed {\bf tree above $r$}, which we call $\tree(r)$, whose nodes are the integers $m\in \Z/L\Z$, such that there exists a directed path from $m$ to $r$ which does not contain any other roots.  We note that $r\in \tree(r)$, and we also call $\tree(0)$ the {\bf mother tree}.
\end{defn}

\begin{thm}[Mother tree]
    \label{thm: The Mother Tree}
$k\in\tree(0)\iff L_n \mid k$. Also, $d(L_n,0)\geq d(k\cdot L_n,0)$, with equality when $(k,n)=1.$

\hfill
    \hyperlink{proof of thm: The Mother Tree}{(Proof)}
    $\square$
\end{thm}

\begin{thm}
\label{thm: distance is invariant over +L_n}
    Suppose $i\in Z(n,L)$ is a root contained in a cycle of size $j$ and let $a$ and $b$ be any two nodes in the same connected component as $i$.  Then
    \[
    a\equiv b\Mod {L_n}\iff d(a,i)\equiv d(b,i)\Mod{j}.
    \]
    \hfill
    \hyperlink{proof of thm: distance is invariant over +L_n}{(Proof)}
    $\square$
\end{thm}

\begin{defn}
The {\bf height} $H(m)$ of a node $m\in Z(n,L)$ is the smallest distance between $m$ and any root of $Z(n, L)$. 
We say that a Zolotarev graph $Z(n,L)$ has  {\bf homogeneous height $h$} if all the leaves have the same height $h$. 
We say $Z(n,L)$ has  {\bf height $h$}, if $h$ is equal to the maximum height of any leaf of the graph.
\end{defn}

\begin{lem}
\label{lem: the max height}
    If $Z(n,L)$ has height $h$, then $h=H(1)$.
    \hfill 
    \hyperlink{proof of lem: the max height}{(Proof)} $\square$
\end{lem}

We can therefore denote the height of $Z(n,L)$ by $H(1)$.

\begin{thm}[Height formula]
    \label{thm: homogeneously height}
    $Z(n,L)$ has homogeneous height $h$ $\iff$ 
    $L_n=\frac{L}{(n,L)^h}$.

\hfill
    \hyperlink{proof of thm: homogeneously height}{(Proof)}
    $\square$
\end{thm}

This formula is reminiscent of our first algorithm to calculate $L_n$, described in Example \ref{recursive formula to L_n}.  We note that the height of the leaves is strongly related to the recursive calculation of $L_n$, as  in Example \ref{recursive formula to L_n}.  
We will revisit these ideas  in Corollary \ref{cor: prune Z(n,L)} below.

%%%%%%%%%%%%%%%%%
\subsection{Graphs isomorphisms}

Let $Z(n,L)$ and $Z(m,K)$ be two Zolotarev graphs. We define $G\subset Z(n,L)$ to be isomorphic to $H\subset Z(m,K)$ if there exists a bijection $\psi:G\rightarrow H$ such that
$\psi$ preserves adjacency between any two nodes of $G$.  That is: 
\begin{equation}
\label{Isomorphic function}
a\longmapsto b \iff \psi(a)\longmapsto \psi(b).
\end{equation}
In this case, we'll use the notation $G\simeq H.$

\begin{thm}
    \label{thm: constructing isomorphism of subgraphs}
    Let $v\mid L$, and consider the subgraph  $G\subset Z(n,L)$ whose nodes are multiples of $v$.
     Then:
    \[
    G\simeq Z\left(n,\frac{L}{v}\right)
    \]
    \hfill
    \hyperlink{proof of thm: constructing isomorphism of subgraphs}{(Proof)}
    $\square$
\end{thm}

\begin{cor}[Mother tree structure]
\label{cor: mother tree isomorphic to Z(n,L)Z(L/L_n)}
    \[
    \tree(0)\simeq Z\left(n,\frac{L}{L_n}\right)
    \]
    \hfill
    \hyperlink{proof of cor: mother tree isomorphic to Z(n,L)Z(L/L_n)}{(Proof)}
    $\square$
\end{cor}

\begin{cor}[Pruning]
    \label{cor: prune Z(n,L)}
  Let $G$ be the subgraph of $Z(n,L)$ whose nodes are not leaves.
Then:
    \[
    G\simeq Z\left(n,\frac{L}{(n,L)}\right)
    \]
    \hfill
    \hyperlink{proof of cor: prune Z(n,L)}{(Proof)}
    $\square$
\end{cor}

The next important corollary follows directly from the latter, and gives us a precise relationship between Zolotarev permutations and Zolotarev graphs.

\begin{cor}
\label{cor: roots of Z(n,L)Z(n,L_n)}
    The subgraph of $Z(n,L)$ consisting of all of its roots is isomorphic to $Z(n,L_n)$.
    
    \hfill
    \hyperlink{proof of cor: roots of Z(n,L)Z(n,L_n)}{(Proof)}
    \hyperlink{annex proof of cor: roots of Z(n,L)Z(n,L_n)}{(Proof 2)} $\square$
\end{cor}

%%%%%%%%%%%%%%%%%%%%%%%%%%%%%%%%%%%%%%%%

\medskip
\subsection{An algorithm for Zolotarev graphs}
The latter three corollaries suggest an efficient algorithm to construct a Zolotarev Graph.

\medskip
\noindent
{\bf Algorithm Z (Construction of Zolotarev graphs)}.
\label{Algorithm:Zolotarev graphs}

\smallskip
 {\bf Step $1$}. \ Create the sequence from the recursive computation of $L_n$, as in Example \ref{recursive formula to L_n}, where the terms are given by $x_0=L$, $x_1=\frac{L}{\rm{gcd}(L,n)}$, $\dots, x_{j+1}=\frac{x_j}{{\rm gcd}(x_j,n)}, \dots$ until we finally get $x_h=L_n$. 
 
 The use of $h$ here is not a coincidence, indeed Corollary 
 \ref{cor: prune Z(n,L)} reveal that the number of steps in the recurrence formula for $L_n$ is exactly equal to the height $h$ of the Zolotarev Graph $Z(n,L)$. 

 \smallskip
{\bf Step $2$}. \ Next, construct the Zolotarev graph $Z(n,L_n)$, which in this case is equal to the Zolotarev permutation $\tau(n,L_n)$.  
 We multiply each node of $Z(n,L_n)$ 
 by $\gcd(n,x_{h-1})$, and create  the leaves connected to the nodes by using Lemma \ref{lem: (n,L) arrows}. 

  \smallskip
{\bf Step $3$}. \
Compute recursively each Zolotarev Graph $Z(n,x_j)$, multiplying every node in $Z(n,x_{j+1})$ by $\gcd(n,x_{j})$, and create  the leaves connected to the nodes by using Lemma \ref{lem: (n,L) arrows}.
\hfill $\square$

\medskip \noindent
We remark that due to Algorithm Z, Theorem \ref{thm: homogeneously height} concerning homogeneous height  becomes  more evident.

%%%%%%%%%%%%%%%%%%%%%%%%%%%%%
\begin{example}
\label{example: computing Z(n,L) using prune}
\rm{
Here we compute the Zolotarev graph $Z(6,60)$, depicted in Figure \ref{Z(6,60)}, by using Algorithm Z, in \eqref{Algorithm:Zolotarev graphs}.
We first calculate $L_n$ recursively: 
\[
x_0=60, \ 
x_1=\frac{60}{6}=10, \ 
x_2=\frac{10}{2}=5=L_n.
\]
Because $n=6 \equiv 1 \Mod 5$, we now 
compute $Z(1,5)$. It is clear that we have $5$ nodes, each forming a cycle of length $1$. 
 Now we multiply all the latter nodes by $\gcd(x_1,n)=\gcd(10, 6)=2$, to obtain the nodes $0$, $2$, $4$, $6$ and $8$ as roots. Note that, by Lemma \ref{lem: (n,L) arrows}, if $r$ is a root, then $r+5\longmapsto r$, thus obtaining $Z(6,10)$.

\begin{figure}[H]
\includegraphics[width=14cm]{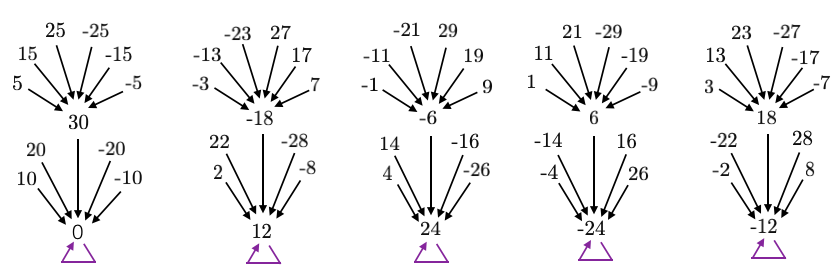}
\caption{The Zolotarev graph $Z(6, 60)$, with its 5 roots, $50$ leaves, $5$ branches, and $5$ connected components.  Its mother tree $\rm{Tree}(0)$ is on the left. } 
\label{Z(6,60)}
\end{figure}

Finally, we multiply every node in $Z(6,10)$ by $\gcd(x_0,n):=\gcd(60, 6)=6$, and connect $6$ new nodes to each leaf in $Z(6,10)$, and $4$ new nodes to each root, obtaining an indegree of $6$ for all the nodes.  For the roots it is easy to compute these nodes because we can obtain them using Lemma \ref{lem: (n,L) arrows}.

\hfill $\square$
}
\end{example}
%%%%%%%%%%%%%%%%%%%%%%%%%%%%%

Example \ref{example: computing Z(n,L) using prune}  strongly suggests the following finer structure, namely that all the trees in a fixed Zolotarev graph are isomorphic to each other.  Indeed this is the case, as the following result shows.

\begin{thm}
\label{thm: trees}
    Given any two roots $r_1, r_2 \in Z(n, L)$, we have
\[
\rm{Tree}(r_1) \simeq	\rm{Tree}(r_2).
\]
\hfill
\hyperlink{proof of thm: trees}{(Proof)}
$\square$
\end{thm}
Theorem \ref{thm: trees} shows us that once we understand the mother tree $\rm{Tree}(0)$, and the cycle structure of the roots, the rest of the graph is completely determined. 

\begin{cor}
\label{cor: number of nodes connected}
    A connected component with a cycle of size $j$, has exactly
$j\cdot\frac{L}{L_n}$ nodes.

\hfill
\hyperlink{proof of cor: number of nodes connected}{(Proof)}
$\square$
\end{cor}

%%%%%%%%%%%%%%%%%%
%%%%%%%%%%%%%%%%%%%%%%%%%%%%%%%%%%%%
%%%%%%%%%%%%%%%%%%

\bigskip
\section{The eigenspaces}
\label{sec:the eigenspaces}

\smallskip
\noindent
Let $r\in Z(n,L)$ be a root in a cycle of size $j$, and let $\omega$ be a primitive $m$'th root of unity, where $m\mid j$.  
We define 
\[
F_{\omega,L,\kappa,n,r}(x)=\sum_{k=0}^{\infty}k^{\kappa-1}a_{\omega,L,n,r}(k)x^k,
\]
where, for all $i\in Z(n,L)$, 
\begin{align*}
    a_{\omega,L,n,r}(i):=
    \begin{cases}
        \omega^{-d(i,r)},&\text{ if $i$ is path-connected to $r$}\\
        0,&\text{ if $i$ is not path-connected to $r$}
    \end{cases}
\end{align*}
\noindent
We note that $a_{\omega,L,n,r}(r)=1$. 
We first enumerate all the cycles whose length is a multiple of $m$:
\begin{equation}
\label{B_m}
    B_m:=\{ \text{ cycle } A \subseteq Z(n, L) \  \mid  \  |A| \equiv 0 \mod m \}.
\end{equation}

To enumerate the eigenfunctions, we may pick the unique smallest element from each of the cycles in $B_m$, as follows.

\begin{thm}[A basis for $E_n(\omega,L,\kappa)$]
\label{thm: basis for En(omega,L,k)}
Let $\omega$ be a primitive $m$'th root of unity.    Then the collection of functions
\begin{equation}
\label{the basis}
   \bigcup_{A \in B_m} 
   \{
   F_{\omega, L, \kappa, n, r} \mid  \  \ r = \min(A) 
   \}
%, \text{ for some } A \in B_m
\end{equation}
is a basis for $E_n(\omega,L,\kappa)$.
\hfill\hyperlink{proof of thm: basis for En(omega,L,k)}{(Proof)} $\square$
\end{thm} 

\begin{example} \label{expanding function in terms of the basis}
\rm{Let 
\[
f(x):= \frac{1}{1-\zeta x} + \frac{1}{1-\zeta^{-1} x},
\]
where $\zeta:= e^{ \frac{2\pi i}{10}}$.  Expanding the latter rational functions into their geometric series: 
\[
f(x) = \sum_{n=0}^\infty \left(\zeta^n + \zeta^{-n} \right)x^n 
=2\sum_{n=0}^\infty 
\cos\left( \frac{2\pi n}{10}\right) x^n.
\]
Applying the operator $U_{9}$, we have:
\begin{align}
U_{9} f(x)&:= U_{9}\left(  \sum_{n=0}^\infty \left(\zeta^n + \zeta^{-n} \right)x^n \right) 
    =  \sum_{n=0}^\infty \left(\zeta^{-n} + \zeta^{(-1)^2 n} \right)x^n = f(x).
\end{align}
Hence $f\in E_9(1,10\cdot m,1)$, for all $m\in\Z_{>0}.$} Choosing $m=3$, we easily see that $\left\{F_{1,30,1,9,r}\right\}$ is a basis for $E_9(1,30,1)$.
\hfill $\square$
\end{example}

%%%%%%%%%%%%%%%%%%

\begin{cor}[Level reduction for eigenspaces]
\label{cor: E_n(L)=E_n(L_n)}
    $E_n(\omega,L,\kappa)=E_n(\omega,L_n,\kappa)$
    \hfill
    \hyperlink{proof of cor: E_n(L)=E_n(L_n)}{(Proof)}
    $\square$
\end{cor}

\begin{thm}[The dimension of $E_n(\omega,L,\kappa)$]
\label{thm:precise dimension of En(omega,L,k)} 
Let $n$ and $L$ be positive integers, $c:= \ord_{L_n}(n)$, $m$ a positive divisor of $c$ and $\omega$ any primitive $m$'th root of unity. Then,
\begin{equation}
   \dim E_n(\omega, L, \kappa) 
   =\sum_{j=1}^{c}b_{jm},
\end{equation}
where $b_j$ is the number of cycles of length $j$ in $Z(n, L)$, and may be computed using \eqref{formula for b_j, Zolotarev graphs}.

\hfill
\hyperlink{proof of thm:precise dimension of En(omega,L,k)}{(Proof)}
\hyperlink{anex proof of thm:precise dimension of En(omega,L,k)}{(Proof 2)}
$\square$
\end{thm}
\begin{example}
\rm{Here we compute the dimension of $E_n(\omega,L,\kappa)$ for $n=5, L=6$, $\omega=-1$, and any weight $\kappa$. Note that 5 is a primitive root modulo 6, so that $c:= \ord_6(5) = 2$. 
Here
\[
\tau(5,6)=(0) \ (1 \ 5) \ (2 \ 4) \ (3),
\]
so we have $b_1=2$ and $b_2=2$. Using $m=2$ (because $\omega=-1$), Theorem \ref{thm:precise dimension of En(omega,L,k)} gives us:
\[
 \dim   E_5(-1, 6, \kappa) 
   =\sum_{j=1}^{1}b_{2j} =2.
\]

On the other hand, proceeding from first principles, suppose that 
$f(x):=\sum_{k=0}^{\infty}a(k)x^k \in E_5(-1,6,\kappa)$.
Then 
    \[\begin{cases}
        a(k)=a(k+6)\\
        a(5k)=-a(k)
    \end{cases}\]
In this case, we have $a(0)=-a(0)=0$, $a(3)=-a(3)=0$ and the following system: 
    \[\begin{cases}
        a(1)=-a(5)\\
        a(2)=-a(4)
    \end{cases},\]
leaving us with two free parameters among these coefficients.   Consequently, $\dim  E_5(-1, 6, \kappa)=2$ again, 
giving us an independent confirmation of Theorem \ref{thm:precise dimension of En(omega,L,k)}.} 
\hfill $\square$
\end{example}
\begin{example}
    \rm{In this example, we compute the dimension of 
    $E_3(e^{\frac{2\pi ia}{6}},7,1)$. Note that 3 is a primitive root modulo 7 (i.e. $c=6$).
    We begin by computing $\tau(3,7):$
    \[(0) \ (1 \ 3 \ 2 \ 6 \ 4 \ 5)\]
    We therefore have $b_1=1$ and $b_6=1$. If $f:=\sum_{k=0}^{\infty}a(k)x^k \in E_3(e^{\frac{2\pi ia}{6}},7,1)$, then 
    \[\begin{cases}
        a(k)=a(k+L)\\
        a(3k)=\omega a(k)
    \end{cases}\]
    In that case, we have $a(0)=\omega a(0)=0$ and 
    \[a(1)=\omega a(5)=\omega^2 a(4)=\omega^3 a(6)=\omega^4 a(2)=\omega^5 a(3)=\omega^6 a(1)\]
As all coefficients $a(k)$ are completely determined by $a(1)$ (i.e., there are no "free variables"), we have
\[\dim E_3(e^{\frac{2\pi ia}{6}},7,1)=\begin{cases}
    1, &\text{if } a\neq 0\\
    2, &\text{if } a=0
\end{cases}
\]}
\hfill $\square$
\end{example}

\begin{example}
    \rm{In this example, we compute the dimension of 
    $E_3(i,7,1)$.
    Again, we have have $b_1=1$ and $b_6=1$. If $f:=\sum_{k=0}^{\infty}a(k)x^k \in E_3(i,7,1)$, then 
    \[\begin{cases}
        a(k)=a(k+L)\\
        a(3k)=i a(k)
    \end{cases}\]
    In that case, we have $a(0)=i a(0)=0$ and 
    \[a(1)=i a(5)=- a(4)=-i a(6)= a(2)=i a(3)=-i a(1)=0\]
Then $f\in E_3(i,7,1)\iff f=0$, therefore 
\[\dim E_3(i,7,1)=0
\]}
\hfill $\square$
\end{example}

\noindent Next, we define the vector space spanned by the functions that are eigenfunctions of $U_n$:
\begin{defn}
Let $n, L$ and $\kappa$ be positive integers. Then,

\[
S_n(L,\kappa):=\left<f\mid U_n(f)
=n^{\kappa-1}\omega f, \text{for some } \omega\in \C
\setminus \{0\},  \text{level}(f)\mid L \text{ and weight}(f)=\kappa\right>,
\]
\end{defn}
 It follows from this definition 
 that: 
\[
S_n(L,\kappa)=\bigoplus_{\omega\in\C
\setminus \{0\}}E_n(\omega,L,\kappa),
 \]

\begin{thm}
\label{thm:precise dimension of Sn(L,k)}
We have:
\begin{equation}
   \dim    S_n(L, \kappa) =L_n,
\end{equation}
where $L_n$ was defined in \eqref{def:L_n}.
\hfill
\hyperlink{proof of thm:precise dimension of Sn(L,k)}{(Proof)}
\hyperlink{anex proof of thm:precise dimension of Sn(L,k)}{(Proof 2)}
$\square$
\end{thm}
%%%%%%%%%%%%%%%%%%

\begin{cor}
\label{cor: S_n(L)=R(L_n)}
    \[
    S_n(L,\kappa)=S_n(L_n,\kappa)=\RR(L_n,\kappa)
    \]
    \hfill
    \hyperlink{proof of cor: S_n(L)=R(L_n)}{(Proof)} $\square$
\end{cor}

\begin{cor}
\label{cor: dim S_n multiplicative and periodic}
   Let $g(n, L):=\dim S_n(L,\kappa)$.  Then $g$ is a periodic function of $n$.  Moreover, $g$ is a completely multiplicative function of $L$.  
   \hfill
   \hyperlink{proof of cor: dim S_n multiplicative and periodic}{(Proof)} $\square$
\end{cor}

%%%%%%%%%%%%%%%%%%%%%%
\section{Simultaneous eigenfunctions}\label{sec: Simultaneous eigenfunctions}

We define the vector space spanned by the simultaneous eigenfunction, namely those rational functions $f\in \RR$ that enjoy $U_n(f)=\lambda f$, for all $n$. Here we allow $\lambda\in\C$, including the case $\lambda=0$.
We study the vector space generated by all of the {\bf simultaneous eigenfunctions}, namely:
\[
V(L,\kappa)=\left<f\mid \text{for all }n, \exists\lambda\in\C \text{, such that } U_n(f)=\lambda f \text{ and level}(f)\mid L\right>.
\]
It follows from the definitions above that
\[
V(L,\kappa)=\left<\bigcap_{n\geq0}\left(\bigcup_{\lambda\in\C}E_n(\lambda,L,\kappa)\right)\right>
\]

We note that there exist functions $f\in V(L,\kappa)$ which are not in $S_n(L,\kappa)$, because $\lambda=0$ is allowed in $V(L,\kappa)$, but not in $S_n(L,\kappa)$.
\medskip

\noindent Given an eigenfunction $f$ of $U_n$, it is natural to ask for which other integers $m$ is it true that $f$ is an eigenfunction of $U_m$. To this end, we introduce the following set.  

\begin{defn} 
\[\mathcal{C}(f):=\left\{n \in \Z_{\geq 1} \mid U_n(f)=\lambda f, \lambda\in\C\right\}\]
\end{defn}
\begin{defn} 
We define $\omega_f:\mathcal{C}(f)\rightarrow\C$ to be the function that satisfies 
\[
U_n(f)=n^{\kappa-1}\omega_f(n) f.
\]
\end{defn}

\bigskip
\noindent
Next, we observe an initial connection of eigenfunctions to complex Dirichlet characters.

\begin{lem} \label{lem: Defining C(f) and lambda_f}
We fix a positive integer $n$, and let $f$ be an eigenfunction of $U_n$, with level $L$ and weight $\kappa$. 
Suppose that $n,m\in\mathcal{C}(f)$. Then,
the following hold:
    \begin{enumerate}[(a)]
        \item $a(k\cdot m)=\omega_f(m)\cdot a(k)$, for all $k\geq0$
        \item $n\cdot m\in\mathcal{C}(f)$
        \item $m+L\in\mathcal{C}(f)$
        \item There exists a Dirichlet character $\chi$ modulo L, such that $\omega_f=\chi|_{\mathcal{C}(f)}$.
    \end{enumerate}
\hfill
\hyperlink{proof of lem: Defining C(f) and lambda_f}{(Proof)}
$\square$
\end{lem}

We also define the finite collection of $\phi(L)$ functions:
\begin{equation}
\label{set of the character series}
\mathfrak{B}_{L,\kappa}:= 
\left\{
f(x):= \sum_{k=0}^\infty 
\chi(k) k^{\kappa} x^k \mid 
\chi \text{ is a Dirichlet character mod } L
\right\}.
\end{equation}

\begin{cor}
\label{cor: generating V_L}
Let $\overline{L}=$level($f$).
    If $f\in V(L,\kappa)$ then $\omega_f=\chi$, a Dirichlet character modulo $\overline{L}$.
    
\hfill
\hyperlink{proof of cor: generating V_L}{(Proof)}
$\square$
\end{cor}

\begin{lem}
\label{(M+kc,L)=1}
Suppose  $L=p\cdot c$, and $M$ is any integer with $(M,L)=1$.  
The following are equivalent:
\begin{enumerate}[(a)]
    \item $p\mid c$.
    \item  For all $k\in \Z_{\geq 0}$, we have $\left(M+k\cdot c,L\right)=1$.
\end{enumerate}  \hfill
\hyperlink{proof of (M+kc,L)=1}{(Proof)}
$\square$
\end{lem}

\noindent
In order to prove 
Theorem \ref{thm:precise dimension of V_L,k} below,
we will require the following technical lemma, which uses some new notions, as follows.  
 Suppose we have the factorization 
    $L=\prod_{i=1}^{n}q_i$, where 
    $q_i:= p_i^{\alpha_i}$, and the $p_i$'s are distinct primes. 
We define:
\begin{equation}
\label{def:A_L}
    A_L:=\left\{
    \prod_{j\in S}q_j
    \mid S\subset \{1, \dots n\}
    \right\}.
\end{equation}

\begin{lem}
    \label{lem: Defining A_l}
    Then, the following hold:
    \begin{enumerate}[(a)]
        \item $a\in A_L\iff a\mid L \text{ and }(a,\frac{L}{a})=1$
        \item $A_L = 
        \{c\; \mid \; L\equiv0\pmod{c} \text{ and } \forall a>1\text{ such that }c\equiv 0 \pmod a\Rightarrow a\cdot c\nmid L\}$
        \item $A_L=(1+q_i)\cdot A_{\frac{L}{q_i}}$
        \item $\left| A_L\right|=2^{n}$
        \item $A_L=\left\{L_n\mid \forall n\in\N\right\}$
    \end{enumerate}
    \hfill
    \hyperlink{proof of lem: Defining A_l}{(Proof)}
    $\square$
\end{lem}
We note that the elements of $A_L$ are called {\bf unitary divisors} of $L$, and have been studied, for example by Cohen \cite{Cohen}.
The following theorem uses the notation in \eqref{set of the character series} and \eqref{def:A_L}.
\begin{thm}[Basis for $V(L, \kappa)$]
\label{thm:precise dimension of V_L,k}
  
  $\dim V(L,\kappa)$ is a multiplicative function of $L$.
  Moreover,  
\begin{equation}
\label{Basis of V(L,kappa)}
\bigcup_{M\in A_L}\mathfrak{B}_{M,\kappa} \text{ is a basis for } V(L,\kappa)
\end{equation} 

and we have the following:
  \begin{enumerate}[(a)]
  \item \label{dimension of V(L,kappa), part a}
\[
\dim V(L,\kappa)=\prod_{i=1}^{n}\left(\phi(p_i^{\alpha_i})+1\right).
\]
\item \label{dimension of V(L,kappa), part b}
\[
\dim V(L,\kappa)=\sum_{M\in A_L} \phi(M).
\]
\item \label{dimension of V(L,kappa), part c}In particular, $\dim V(L,\kappa)=L \iff  L$ is square free.
\end{enumerate}
\hfill
\hyperlink{proof of thm:precise dimension of V_L,k}{(Proof)}
$\square$
\end{thm}

We remark that 
rational functions with completely multiplicative coefficients are contained inside the vector space $V(L,\kappa)$ of simultaneous eigenspaces.

%%%%%%%%%%%%%%%%

%%%%%%%%%%%%%%%%%%%%%%%%%%%%%%%%%%%%%%%%%%%%%%%%%%%%%%%%

\begin{comment}
    What is the dimension of the "old" $V_L^{simult}$ how we defined in the start, it means, with the Level$(f)=L$ and not a divisor.

    Maybe smaller that $\phi(L)$?

    If $L$ is prime it must be easy to solve, but if not, must be difficult and have relation with the primitivity of Dirichlet Characters, a much more hardy question. Maybe we can solve a little using our formulas, and also apply it in some thing.
\end{comment}

%%%%%%%%%%%%%%%%%%%%%%%%%%%%%%%%%%%%%%%%%%%%%%%%%%%%%%
%%%%%%%%%%%%%%%%%%%%%%%%%%%%%%%%%%%%%%%%%%%%%%%%%%%%%%

\section{The Artin conjecture}\label{sec: The Artin conjecture}
In this brief section we give an equivalence between the Artin conjecture \cite{Gupta.Murty}, and certain Eigenspaces that appeared in Theorem \ref{thm:precise dimension of En(omega,L,k)}.
\begin{conj}[Artin, 1927]
Suppose $n$ is an integer that is not a perfect square, and $n\neq-1$.  Then $n$ is a primitive root mod $p$ for infinitely many primes $p$.
\end{conj}
\noindent 
If the conclusion above is true, we say that {\bf the Artin conjecture is true for $n$}. Gupta and Murty \cite{Gupta.Murty} proved in $1983$ that unconditionally, there is a set of $13$ integers such that the Artin conjecture is true for at least one of these integers.  Their proof also showed that the Artin conjecture is true for almost all $n$. Indeed one of their ideas was to use sieve theory to produce primes $p$ such
that $p-1$ has only a few prime factors. 
In $1985$ Heath-Brown \cite{Heath-Brown} reduced the size of this set to $3$ integers, thereby proving that there are at most $2$ primes for which the Artin conjecture fails.  It is not yet known, unfortunately, how to check whether any given prime is one of these two possible exceptions.  In \cite{Jensen.Murty}, Jensen and Murty studied an analogue of Artin's conjecture for polynomials over finite fields, and simultaneously simplified some of the arguments in the literature.  
We can now rephrase the Artin conjecture in terms of $1$-dimensional eigenspaces of $U_n$.
\begin{thm}[An equivalence for the Artin conjecture]
The following are equivalent:
\begin{enumerate}[(a)]
\item \label{part a of Artin}
Fix any positive integer $\kappa$, and $\omega:= e^{ \frac{2\pi i}{p-1}}$.
Then there are infinitely many primes $p$ such that 
\[
\dim E_n(\omega,p,\kappa)=1.
\]
\item  
\label{part b of Artin}
\begin{align}
\frac{1}{p-1}
\sum_{d \mid p-1} 
\mu\left(\frac{p-1}{d}\right)
\gcd(n^d-1, p)=1,
\end{align}
for infinitely many primes $p$.
\item 
\label{part c of Artin}
The Artin conjecture is true for $n$.
\end{enumerate}
\end{thm}
\begin{proof}
By definition $n$ is a primitive root mod $p \iff \ord_p(n) = p-1$. Using our definition of $b_j$ as the number of cycles of length $j$ in the permutation $\tau(n, L)$, we have 
$b_{p-1}=1$ while $b_j=0$ for all 
$1< j\not=p-1$. 
By Theorem \ref{thm:precise dimension of En(omega,L,k)}, 
$\dim E_n(\omega,p,\kappa)
=\sum_{j=1}^{\frac{c}{m}}
b_{jm}=1$, proving the equivalence of  \eqref{part a of Artin} and \eqref{part c of Artin}.
The equivalence of \eqref{part b of Artin} 
and \eqref{part a of Artin}
follows from  equation 
\eqref{formula for b_j, Zolotarev graphs}, which is an explicit formula for all $b_j$.
\end{proof}

%%%%%%%%%%%%%%%%%%%%%%%%%%%
\bigskip
\section{The kernel}\label{sec: The Kernel}

Given any positive integer $n$, we classify all rational functions $f\in \RR(L,\kappa)$ that belong to the kernel of $U_n:\RR(L,\kappa) \rightarrow \RR(L,\kappa)$.  That is, we define:
\[
\ker(U_n) := \{ f\in \RR(L,\kappa) \mid U_n f(x) = 0, \text{ for all } x \text{ near } 0\}.
\]
We note that if we were to consider 
$U_n:\RR \rightarrow \RR$, then 
the structure theorems that we proved, concerning eigenfunctions and eigenvalues of $U_n$, no longer hold. In particular, the following example shows that there are functions in the kernel of $U_n:\RR \rightarrow \RR$ that do not belong to $\RR(L,\kappa)$, for any $L$ and $\kappa$.

\begin{example}
\label{example: weight does not exist}
\rm{The following function does not lie in $\RR(L,\kappa)$, so this example only serves to show what can happen when a function is in the kernel of $U_n$, but lies in $\RR\setminus\RR(L,\kappa)$. 

Let $f\in \RR$ 
be the rational function defined by:
\[
f(x):= \sum_{k=1}^\infty 
\left(
1_{2 \, {\rm mod} 9}(k) \cdot k  +  
1_{7 \, {\rm mod} 9}(k)
\right)x^k,
\]
where we use the notation of \eqref{Indicator Function}.
We note that $a_0=a_3=a_6=0$, and it is easily verified that  $U_6(f)=0$.
We emphasize that for functions in the kernel, the notion of ``weight'' does not make sense, as we see in this example.} 
\hfill $\square$
\end{example}
The latter example shows that it makes sense to restrict the action of $U_n$ to  $\RR(L,\kappa)$, where the notions of level and weight of its eigenfunctions are well defined.

\begin{lem}[The kernel in terms of its Zolotarev graph]
\label{lem: equivalence for the kernel}
Let $U_n: \RR(L,\kappa) \rightarrow \RR(L,\kappa)$ and let 
$f\in\RR(L,\kappa)$. Then 
\[
f\in \ker(U_n) \iff a(k)=0,  \ \text{for all branches and roots } k\in Z(n,L).
\]
\hfill
\hyperlink{proof of lem: equivalence for the kernel}{(Proof)}
 $\square$
 \end{lem}

\begin{example}
\rm{
Let $f\in \RR(4, 1)$. In other words, the
Taylor series coefficients $a(k)$ of $f$ have period $L=4$:
$a(0)=0, \ a(1)=3,\    a(2)=0, \ a(3)=17$.
With $n=2$, we apply $U_n$, and it's easy to see that $U_2(f)(x) = 0$.
Figure \ref{L=4,n=2} illustrates this Zolotarev graph $Z(2, 4)$.

\begin{figure}[ht]
\includegraphics[width=2cm]{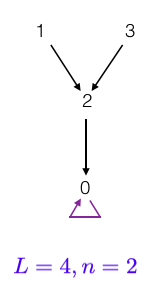}
\caption{The Zolotarev graph $Z(2, 4)$ has 
$2$ leaves, $1$ branch, and $1$ root. The nodes
 $1$ and $3$ are leaves in $Z(2,4)$, which means that the Taylor series coefficients $a(1)$ and $a(3)$ represent free variables. On the other hand, a sufficient and necessary condition for $f\in \ker(U_2)$ is that $a(0)=a(2)=0$ (the roots and branches in $Z(2,4)$).  Therefore, it is clear that $\dim\ker(U_2)=2$.} 
\label{L=4,n=2}
\end{figure}

}
\hfill $\square$
\end{example}

\noindent

\begin{cor}[Dimension of the kernel]
    \label{cor:precise dimension of ker(U_n)}
   Let $U_n:\RR(L,\kappa)\rightarrow\RR(L,\kappa).$ Then,
    \[
    \dim \ker(U_n)=L - \frac{L}{(n,L)}.
    \]
\hfill
\hyperlink{proof of cor:precise dimension of ker(U_n)}{(Proof)}
$\square$
\end{cor}

%%%%%%%%%%%%%%%%%%%%%%%%%%%%%%%%%%%%%%%%%%%%%
%  The following result is now a main result of the paper,
%  giving equivalances for diagonalizability
%%%%%%%%%%%%%%%%%%%%%%%%%%%%%%%%%%%%%%%%%
The following result gives equivalences for diagonalizability.

\begin{thm}[Diagonalizability of $U_n$]
\label{thm: diagonalizable of U_n}
The following conditions are equivalent:
\begin{enumerate}[(a)]
\item
$U_n:\RR(L,\kappa)\rightarrow\RR(L,\kappa)$ is diagonalizable.
\item   $\frac{L}{(L,n)}=L_n$
\item $Z(n,L)$ has no branches.
\item $n$ is a root in $Z(n,L)$
\item $H(1)=1$
\end{enumerate}
\hfill
\hyperlink{proof of thm: diagonalizable of U_n}{(Proof)}
$\square$
\end{thm}

%%%%%%%%%%%%%%%%%%%%%%%%%%%
%%%%%%%%%%%%%%%%%%%%%%%%%%%%%%%%%%%
%%%%%%%%%%%%%%%%%%%%%%%%%%%%%%%%%%%

\section{Proofs for section \ref{sec:intro} Introduction}

\bigskip
\subsection{Proof of Theorem \ref{thm:first result} and Lemma \ref{gcd(n, L)=1 for an eigenfunction}}

\begin{proof}
\hypertarget{proof of thm:first result} (of Theorem \ref{thm:first result})
By \cite{Stanley}, Chapter 4, we may write the Taylor series coefficients of any function $f\in\RR$ as 
\[
a(k)=\sum_{j=1}^{d}P_j(k)r_j^k, 
\]
where degree $(P_j(k))= m_j-1$, $m_j$ is the multiplicity of each distinct pole $\frac{1}{r_j}$ of $f$ and $d$ is the number of poles of $f$. Therefore the Taylor series coefficients of $U_n(f)$
are
\[
a(kn)=
\sum_{j=1}^{d}P_j(kn)r_j^{kn}
=
\sum_{j=1}^d P_j(kn)(r_j^n)^k,
\]
which shows that the poles of $U_n(f)$ are precisely 
$\left\{ \frac{1}{r_1^n}, \frac{1}{r_2^n}, \dots, \frac{1}{r_d^n}\right\}$.
On the other hand, our assumption that 
$U_n(f) = \lambda f$ implies that the rational functions $f$ and 
 $U_n(f)$ have the same set of poles, namely $\left\{ \frac{1}{r_1}, \dots, \frac{1}{r_d}\right\}$.  We now have:
\begin{equation}
    \left\{\frac{1}{r_1}, \frac{1}{r_2}, \cdots, \frac{1}{r_d}\right\}=
    \left\{\frac{1}{r_1^n},\frac{1}{r_2^n}, \cdots, \frac{1}{r_d^n}\right\},
\end{equation}
implying that the latter two sets must be permutations of each other.  Considering the orbit of the pole $\frac{1}{r_i}$ under $U_n$, we know that the subset of poles of $f$ defined  
by $\left\{ \frac{1}{r_i}, \frac{1}{r_i^n}, \frac{1}{r_i^{n^2}}, \dots \right\}$ must be a finite set, so that for some integers $j, k$, we have 
$\frac{1}{r_i}^{n^j}=\frac{1}{r_i}^{n^k}$.  Therefore $\frac{1}{r_i}^{n^j-n^k}=1$, so that $\frac{1}{r_i}$ is a root of unity and then also $r_i$ is a root of unity.
\end{proof}

\begin{proof}
\hypertarget{proof of gcd(n, L)=1 for an eigenfunction} (of Lemma \ref{gcd(n, L)=1 for an eigenfunction})
The Taylor series coefficients of any function $f\in\RR$ are given by \cite{Stanley}
\[
a(k)=\sum_{j=1}^{d}P_j(k)r_j^k, 
\]
where degree$(P_j(k))= m_j-1$, $m_j$ is the multiplicity of each distinct pole $\frac{1}{r_j}$ of $f$ and $d$ is the number of poles of $f$.  
If $f$ is an eigenfunction of $U_n$, with level $L$, then we know by Theorem \ref{thm:first result} that all of its poles $\frac{1}{r_j}$, (and thus $r_j$) are $L$'th roots of unity. Indeed, there is a smallest integer $L$ with this property.   Therefore:
\[
a(k)=\sum_{j=1}^{d}P_i(k)
e^{\frac{2\pi i l_j k }{L} },
\]
for some integers $l_j$.  

Now we suppose to the contrary that $\gcd(n, L)=M > 1$, so we have $n = n_0\cdot M, L = L_0\cdot M$, for some integers $ n_0 < n, L_0 < L$. 
The action of $U_n$ on $f$ gives us the relation:
\begin{align}
a(k) \lambda = a(kn) = 
\sum_{j=1}^{d}P_i(nk)
e^{\frac{2\pi i l_j k n }{L} }
=
\sum_{j=1}^{d}P_i(k)
e^{\frac{2\pi i l_j k n_0 }{L_0} },
\end{align}
implying that the poles of $f$ are all $L_0$´roots of unity, with $L_0 < L$, a contradiction.
\end{proof}

%%%%%%%%%%%%%%%%%%%%%%%%%%%%%%%%%%%
%%%%%%%%%%%%%%%%%%%%%%%%%%%%%%%%%%%
\medskip
\subsection{Proof of
%Lemma \ref{lemma:How to write the Taylor Series of an eigenfunction}, of 
Theorem
\ref{thm:second result}, 
Theorem \ref{thm:eigenvalues and kappa}
and Corollary \ref{corollary:How to write an eigenfunction}}
\hypertarget{proof of thm:second result}.

\noindent Here, we'll prove
%Lemma \ref{lemma:How to write the Taylor Series of an eigenfunction}, 
Theorem \ref{thm:second result}, 
Theorem \ref{thm:eigenvalues and kappa} and
Corollary \ref{corollary:How to write an eigenfunction}
simultaneously.

\begin{proof}
    We may write the Taylor series coefficients of any function $f\in\RR$ as 
\[
b(k)=\sum_{j=1}^{d}P_j(k)r_j^k, 
\]
where degree$\left(P_j(k)\right)= m_j-1$ and $m_j$ is the multiplicity of each distinct pole $\frac{1}{r_j}$ of $f$.
By assumption, we have 
\[
U_n(f)=\sum_{k=0}^{\infty}b(nk)x^k=\lambda \sum_{k=0}^{\infty} b(k)x^k.
\]
Then, for all $k\in\Z_{\geq0}$,
\begin{equation}\label{a(nk)=lambda a(k)}
    b(nk)=\lambda b(k)
\end{equation}
Let $L=\rm{level}(f)$. By Theorem \ref{thm:first result}, for all $i=1,2,\dots,d$, we have $r_i^L=1$. Also by Lemma \ref{gcd(n, L)=1 for an eigenfunction}, $(n,L)=1$, hence we can define $c:=\ord_L(n)$.
Now, iterating $U_n$ $c$ times to $f$, we have
\[
b(kn^{c})=
\sum_{j=1}^{d}P_j(kn^{c})r_j^{kn^{c}}
=
\lambda^{c} \sum_{j=1}^d P_j(k)(r_j)^k.
\]
On the other hand, $n^{c}\equiv 1 \mod L$. Therefore, we have
\[
\sum_{j=1}^{d}P_j(kn^{c})\left(r_j^{n^{c}}\right)^k
=
\sum_{j=1}^{d}P_j(kn^{c})r_j^k
=
\lambda^{c} \sum_{j=1}^d P_j(k)(r_j)^k.
\]
By uniqueness of the finite Fourier series, for all $j=1,2,\dots, d$
\begin{equation}
\label{eq: Pj(kn^c)}
P_j(kn^{c})=\lambda^{c} P_j(k).
\end{equation}
Now, because \eqref{eq: Pj(kn^c)} holds for every integer $k$, we have that the coefficients of both polynomials are indeed equal. Let $P_j(k)=\sum_{i=0}^{m_j-1}a_{i,j}k^i.$ The latter observation can be expressed by the following system:
\[
\begin{cases}
    a_{0,j}=\lambda^{c}a_{0,j}\\
    a_{1,j} n^{c}=\lambda^{c}a_{1,j}\\
    \vdots\\
    a_{i,j}n^{ic}=\lambda^{c}a_{i,j}\\
    \vdots\\
    a_{m_j-1,j}n^{(m_j-1)c}=\lambda^{c}a_{m_j-1,j}
\end{cases}
\]

Note that $a_{m_j-1,j}\neq0$, because degree$(P_j(k))=m_j-1$.  Thus, $n^{(m_j-1)c}=\lambda^{c}$

But, because $n>1$, $\lambda^{c}\neq n^{ic}$, for all $i\neq m_j-1$, and, therefore, $a_{i,j}=0$. Indeed, 
\begin{equation}
    \label{Series expansion}b(k)=\sum_{j=1}^{d}P_j(k)r_j^k=\sum_{j=1}^{d}C_jk^{m_j-1}r_j^{k},
\end{equation}
where $C_j=a_{m_j-1,j}$.

Finally, from the equality $n^{(m_j-1)c}=\lambda^{c}$, which holds for all $j=1,2,\dots,d$, and from the fact that $\lambda$ is a constant independent of $j$, we conclude that $\kappa:=m_1=m_2=\dots=m_d$, proving Theorem \ref{thm:second result}.
Moreover, for some $c$'th root of unity $\omega$,
\[
\lambda=\omega n^{\kappa-1},
\]
which concludes the proof of Theorem \ref{thm:eigenvalues and kappa}.

Now, returning to the series expansion \eqref{Series expansion}, we replace all the multiplicities $m_j$ by $\kappa$, getting
\[
b(k)=\sum_{j=0}^{d}k^{\kappa-1}C_jr_j^k=k^{\kappa-1}\sum_{j=0}^{d}C_jr_j^k.
\]
Let $a(k)=\frac{b(k)}{k^{\kappa-1}}$. It follows from Theorem \ref{thm:first result} that $a(k+L)=a(k)$, because
\begin{equation}
    a(k+L)=\sum_{j=0}^{d}C_jr_j^{k+L}=\sum_{j=0}^{d}C_jr_j^kr_j^L=\sum_{j=0}^{d}C_jr_j^k=a(k).
\end{equation}
Moreover, suppose there exists an integer $M<L$ such that $a(k+M)=a(k)$ for all $k$, i.e.
\[
a(k+M)=\sum_{j=0}^{d}C_jr_j^{k+M}=\sum_{j=0}^{d}\left(C_jr_j^{M}\right)r_j^{k}=\sum_{j=0}^{d}C_jr_j^{k}=a(k)
\]
        Then, by the uniqueness of the finite Fourier series, all the coefficients must be equal. In particular, for $k=0$,
        \[
        C_jr_j^M=C_j\Rightarrow r_j^M=1\,\forall j,
        \] a contradiction, as this would imply that the level of $f$ is $M\not=L$, and therefore concluding the proof of Corollary \ref{corollary:How to write an eigenfunction}.
\end{proof}

%%%%%%%%%%%%%%%%%%%%%%%%%%%%%%%%%%%
%%%%%%%%%%%%%%%%%%%%%%%%%%%%%%%%%%%
\medskip
\subsection{Proofs of 
%Corollary \ref{corollary:How to write an eigenfunction} and 
Corollary
\ref{cor:the spectrum of U_n} and Corollary \ref{cor: image of phi function}
}

%\begin{proof}
%\hypertarget{proof of Corollary:How to write an eigenfunction}(of Corollary \ref{corollary:How to write an eigenfunction})
%       Using the structure of the eigenfunctions from Lemma \ref{lemma:How to write the Taylor Series of an eigenfunction} , we write $f$ as \[f(x)=\sum_{k=0}^{\infty}b(k)x^k\] where 
%        \[
%        b(k)=\sum_{j=0}^{d}k^{m_j-1}C_jr_j^k
%        \]
%        But, from \ref{thm:second result}, $m_j=\kappa$, for all $j$, and $r_j$ are $L$´th roots of unity, where $L$ is the level of $f$. Therefore:
%        \[
%        b(k)=k^{\kappa-1}\sum_{j=0}^{d}C_jr_j^k
%        \]
%        Define $a(k)=\frac{b(k)}{k^{\kappa-1}}$. It follows from Theorem \ref{thm:first result} that $a(k+L)=a(k)$, because
%        \begin{equation}
%            a(k+L)=\sum_{j=0}^{d}C_jr_j^{k+L}=\sum_{j=0}^{d}C_jr_j^kr_j^L=\sum_{j=0}^{d}C_jr_j^k=a(k).
%        \end{equation}
%        Now, suppose there exists an integer $M<L$ such that $a(k+M)=a(k)$ for all $k$, i.e.
%\[
%a(k+M)=\sum_{j=0}^{d}C_jr_j^{k+M}=\sum_{j=0}^{d}\left(C_jr_j^{M}\right)r_j^{k}=\sum_{j=0}^{d}C_jr_j^{k}=a(k)
%\]
%        Then, by the uniqueness of the finite Fourier series, all the coefficients must be equal. In particular, for $k=0$,
%        \[
%        C_jr_j^M=C_j\Rightarrow r_j^M=1\,\forall j,
%        \] an absurd, as this would imply that the level of $f$ is $M\not=L$.
%    \end{proof}
    
\begin{proof}
    \hypertarget{proof of cor:the spectrum of U_n}(of Corollary \ref{cor:the spectrum of U_n})
We define:
    \begin{equation}
\mathfrak{A}_n:=  
\{\lambda(m,L,\kappa):=    
e^{ \frac{2\pi i m}{c}} n^{\kappa-1}  \mid
m \in \Z, \ \kappa, L \text{ are  positive integers},
 \text{ and } (n,L)=1
\},
\end{equation}
where $c:=\ord_L(n)$.
Then by Theorem \ref{thm:eigenvalues and kappa} and Lemma \ref{gcd(n, L)=1 for an eigenfunction}, we know that $\rm{Spec}\left(U_n\right)\subseteq\mathfrak{A}_n$.

To prove that $\mathfrak{A}_n \subseteq \rm{Spec}\left(U_n\right)$, it remains to show that for each $m \in \Z$, and $\kappa, L$
positive integers with $(n,L)=1$, there exists a function $f$ that satisfies $U_n(f)=\lambda(m,L,\kappa) f$. Fix $\rm{level}(f)=L$ and $\rm{weight}(f)=\kappa$. Then, by Corollary \ref{corollary:How to write an eigenfunction}, we have
    \[
    f(x)=\sum_{k=0}^{\infty}k^{\kappa-1}a(k)x^k.
    \]
    We fix $a(1)=1$ and $\omega:=e^{\frac{2\pi i m}{c}}$.
We define:
\[
    a(k):=
\begin{cases}
    \omega^{v}, &\text{ if  there exists an integer } v\Mod c
      \text{ such that } k \equiv n^{v} \Mod L,  \\
    0, &\text{ otherwise.}
\end{cases}
\]
We first show that $a(k)$ is well defined.  Suppose we have two integers $v_1, v_2$ that satisfy the latter property. If
$v_1\equiv v_2\Mod c$, then $a(k)=\omega^{v_1}=\omega^{v_2}$, because $\omega$ is a $c$'th root of unity. 
Now we apply $U_n$ to $f$:
\[
    U_n(f(x))=\sum_{k=0}^{\infty}(nk)^{\kappa-1}a(nk)x^k.
\]
If $k\equiv n^{v} \Mod L $, then $ nk\equiv n^{v+1} \Mod L$, thus 
 $a(nk)=\omega^{v+1}=\omega a(k)$.  On the other hand, if such an integer $v$ does not exists for $k$, then it does not exist for $nk$.  Therefore $a(nk)=0=\omega a(k)$. Therefore,
    \[
    U_n(f(x))=\sum_{k=0}^{\infty}n^{\kappa-1}k^{\kappa-1}\omega a(k)x^k=\omega n^{\kappa-1}f(x)
    \]
Thus, for each triplet $(m,L,\kappa)$, there exists an eigenfunction $f$ of $U_n$, with eigenvalue $\lambda(m,L,\kappa)$. 
Therefore, we conclude that 
$\mathfrak{A}_n  \subseteq\rm{Spec}
\left(U_n\right)$.
\end{proof}

\begin{proof}
(of Corollary \ref{cor: image of phi function})
\hypertarget{proof of cor: image of phi function}
Let $x\in\Z_{>0}$ be such that $\phi(x)=N$. By Dirichlet´s theorem on primes in arithmetic progressions, there exists a prime $p$ such that $p\equiv1\pmod{N}$, i.e., there exists $j\in\Z_{>0}$ such that $N\cdot j=p-1=\phi(p).$ Let $g$ be a primitive root modulo $p$. We now have $\ord_p(g^j)=N$. Therefore, by Corollary \ref{cor:the spectrum of U_n}, $e^{\frac{2\pi i}{N}}\in\text{Spec}(U_n)$.
\end{proof}

\subsection{Proof of Lemma \ref{lem: precise dimension of R_L,k}}

\begin{proof}
(of Lemma \ref{lem: precise dimension of R_L,k})
\hypertarget{proof of lem: precise dimension of R_L,k} We will prove that $\RR(L,\kappa)$ is a finite subspace of the infinite dimensional vector space $\RR$. It is clear that $\RR(L,\kappa)\subset\RR$. Note also that the null function $0(x)=0, \forall x\in\C$ lies in $\RR(L,\kappa)$, for all $L$ and $\kappa$. Indeed,
$0=\sum_{k=0}^{\infty}k^{\kappa-1}a(k)x^k$, where $a(k)=0$ for every $k$, thus $0\in\RR(L,\kappa)$.
Now, we prove that $\RR(L,\kappa)$ is closed over sum and multiplication by scalar:
For every $f,g\in\RR(L,\kappa)$, we have
\[
(f+g)(x)=\sum_{k=0}^{\infty}k^{\kappa-1}a(k)x^k+\sum_{k=0}^{\infty}k^{\kappa-1}b(k)x^k=\sum_{k=0}^{\infty}k^{\kappa-1}(a+b)(k)x^k.
\]
Note that $(a+b)(k+L)=a(k+L)+b(k+L)=a(k)+b(k)=(a+b)(k)$, i.e., $f+g\in\RR(L,\kappa)$. Also, for every $\alpha\in\C$, we have
\[
\alpha\cdot f=\sum_{k=0}^{\infty}k^{\kappa-1}(\alpha\cdot a)(k)x^k.
\]
Note that $(\alpha\cdot a)(k+L)=\alpha\cdot a(k+L)=(\alpha\cdot a)(k)$, thus $\alpha f\in\RR(L,\kappa)$, proving that $\RR(L,\kappa)$ is a subspace of $\RR$.

It is easy to see that $\dim\RR(L,\kappa)=L$. Indeed, as defined in \eqref{Indicator Function}, the indicator function $1_{a \text{ mod }L}(k)$ of the congruence class $\{n\in\Z\mid n\equiv a\mod L\}$ is the natural basis for the space of all periodic functions with period $L$. In particular, by Corollary \ref{corollary:How to write an eigenfunction}, the Taylor series coefficients $a(k)$ have period $L$. Thus, we have that 
\[
\left\{\sum_{k=0}^{\infty}k^{\kappa-1}1_{a \text{ mod }L}(k)x^k=\sum_{k=0}^{\infty}(a+kL)^{\kappa-1}x^{(a+kL)} \mid 0\leq a<L\right\}
\]
is a basis of $\RR(L,\kappa)$, which proves that $\dim\RR(L,\kappa)=L$.
\end{proof}
%%%%%%%%%%%%%%%%%%%%%%%%%%%%%%%%%%%
%%%%%%%%%%%%%%%%%%%%%%%%%%%%%%%%%%%
\medskip
\subsection{Proofs of Lemma \ref{formula for b_j} and Lemma \ref{lem:defining b_j}}

\begin{proof}
\hypertarget{proof of Lemma: formula for the b_j}(of Lemma \ref{formula for b_j})
For clarity, we first prove the formula for the number of fixed points of $\tau(n, L)$, namely we must show that 
$b_1= \gcd(n-1, L)$, as predicted by formula \eqref{Formula for b_j}. 
Suppose
 that some $k\mod L$ forms a cycle of length $1$, in the permutation $\tau(n, L)$.  This means that 
$nk\equiv k \mod L$.  The latter congruence is equivalent to $L \mid k(n-1)$, which in turn is equivalent to 
$\frac{L}{\gcd(L, n-1)} \mid k$.   Now it's easy to count how many integers $0\leq k < L$ are multiples of 
$\frac{L}{\gcd(L, n-1)}$: there are exactly $\gcd(L, n-1)$ of them, including $0$.  This proves formula \eqref{Formula for b_j} for $j=1$.

Now we extend the argument above, to prove the required formula 
\eqref{Formula for b_j} for all $b_j$.  Suppose we have a cycle of length less than or equal to $j$:  
$(k, \overline{nk}, \dots, \overline{n^{j-1}k})$, which means that 
\[
n^j k \equiv k \mod L,
\]
and it is elementary that the latter congruence is equivalent to 
\[
\frac{L}{\gcd(n^j-1, L)} \mid k. 
\]
The same divisibility criterion also holds for each of the elements in the same cycle.  Therefore, counting the total number of cycles of length $\leq j$, we have:
\[
\sum_{d\mid j} d b_d =
\gcd(n^j-1, L).
\]
Applying M\"obius inversion, we arrive at:
\begin{align*}
j b_j 
&=  \sum_{d\mid j} 
\mu\left(\frac{j}{d}\right)
\gcd(n^d-1, L).
\end{align*}
\end{proof}

\begin{proof}
\hypertarget{proof of lem:defining b_j}
(of Lemma \ref{lem:defining b_j})
\begin{enumerate}[(a)]
    \item Suppose there exists  $k \nmid c={\rm ord}_L(n)$ such that $k$ is the size of a cycle that $n$ creates in $\Z/L\Z$, $i.e.$, there exists an integer $m$ such that 
\[
m\cdot n^k\equiv m \pmod L.
\]
Multiplying the two sides by $n^{c-k}$, we have
\[
m\cdot n^{c}\equiv m\cdot n^{c-k} \pmod L.
\]
Now, we use $n^c\equiv 1 \pmod L$, getting
\[
m\equiv m\cdot n^{c-k} \pmod L.
\] 
Then, by definition, $k$ is the smallest integer $a$ such that $m\equiv m\cdot n^{a} \pmod L$. Thus, $c-k>k$, i.e., $k<\frac{c}{2}$.

Now we can take $c-k\equiv b \pmod k$, where $0<b<k$, since $k\nmid c$.
But, therefore
\[
m\equiv m\cdot n^{b} \pmod L,
\]
a contradiction. \newline
\item It's clear that if $A\cap B=\varnothing$,  then $A\not=B$, since $A\not=\varnothing$.

Now we prove the converse. Let $A=\left\{a\cdot n^i\in\Z/L\Z\mid i\in\Z\right\}$ and $B=\left\{b\cdot n^i\in\Z/L\Z\mid i\in\Z\right\}$. Suppose there exist $i$ and $j$ such that $a\cdot n^i\equiv b\cdot n^j \pmod L$, i.e., $A\cap B\not=\varnothing$. We will show that the sets A and B are, in fact, equal.
 
 For all $m\in\Z$, we have $a\cdot n^{i+m}\equiv b\cdot n^{j+m}\pmod L$. We therefore create the cycle $A=\left\{a\cdot n^{i+m}\in\Z/L\Z\mid  m\in\Z\right\}=\left\{b\cdot n^{j+m}\in\Z/L\Z\mid  m\in\Z\right\}=B$.\newline
 \item The latter proof not only shows that all  different cycles are disjoint, but also proves that if $b_j$ is the number of different cycles of size $j$, then $\sum_{j=1}^{c}jb_j=L$, concluding the proof of the Lemma.
\end{enumerate}

\end{proof}

%%%%%%%%%%%%%%%%%%%%%%%%%%%%%%%%%%%%%%%
%%%%%%%%%%%%%%%%%%%%%%%%%%%%%%%%%%%%%%%

\medskip
\section{Proofs for section \ref{sec: Zolotarev graphs},  Zolotarev Graphs}

\subsection{Proof of Lemma \ref{lem: (n,L) arrows}}

\begin{proof}
\hypertarget{proof of lem: (n,L) arrows}
    (of Lemma \ref{lem: (n,L) arrows})
  Let $c$ be one of the $(n,L)$ solutions to $n\cdot c \equiv k \mod L$. Writing any of these solutions 
 as $c+m$, we have:
    $n\cdot (c+m)=n\cdot c +n\cdot m \equiv k+n\cdot m \mod L$. 
    But,
    \[
    k+n\cdot m\equiv k \Mod L \iff n\cdot m\equiv 0\Mod L\iff \exists b\in\Z \text{ such that }n\cdot m=b\cdot L.
    \]
Letting $n=d\cdot (n,L)$, with $(d,L)=1$, we may therefore write: 
\[
    m=b\cdot\frac{L}{d(n,L)}=\frac{b}{d}\cdot\frac{L}{(n,L)}.
\]
Using $(d,L)=1$, we have 
    \[
    \frac{b}{d}\cdot\frac{L}{(n,L)}\in\Z\implies\frac{b}{d}\in\Z.
    \]
    Thus, $m$ is another solution to $n\cdot c \equiv k \pmod{L}$ if, and only if, $m=c+j\frac{L}{(n,L)}$. Therefore, all of the  $\gcd(n,L)$ solutions are congruent to each other modulo $\frac{L}{(n,L)}.$
\end{proof}

\subsection{Proofs of Theorem \ref{thm: zolotarev structure}, Corollary \ref{cor: number of roots, leaves and branches} and Lemma \ref{lem: m+kL_n is a root?}}

\begin{proof}
\hypertarget{proof of thm: zolotarev structure} (of Theorem \ref{thm: zolotarev structure})
\begin{enumerate}[(a)]
\item Let $m\in Z(n,L)$ and suppose $(n,L)\mid m$. Then, there exists an integer $k$ such that
\newline$m\equiv n\cdot k\Mod L$. Thus, by definition $m$ is not a leaf.
The converse is analogous.

\item Suppose $\frac{L}{L_n}\nmid m$. Note that $\frac{L}{L_n}$ is exactly the product of the powers of primes that divide $L$ and appear in the factorization of $n$ (with eventually different exponents). Thus, there exists an integer $k$´ such that $\frac{L}{L_n}\mid m\cdot n^{k'}$. Then, by definition, $m$ is not a root.

Conversely, suppose $m=k\cdot\frac{L}{L_n}$, for some integer $k$. Indeed, there exists an integer $j$ such that
\begin{equation}
    n\cdot m=n\cdot k\cdot \frac{L}{L_n}\equiv j\cdot\frac{L}{L_n}\Mod{L}.
\end{equation}
Dividing by $\frac{L}{L_n}$, 
\begin{equation}
    n\cdot k\equiv j\Mod{L_n}.
\end{equation}

 Now, because $(n,L_n)=1$, there exists a positive exponent $l$ such that 
 \[n^l\cdot k\equiv k\Mod{L_n}\implies n^l\cdot m\equiv m\Mod{L},\]
 i.e., $m$ is a root in $Z(n,L)$, concluding the proof of (b).

\item follows trivially from the latter two items and the definition of branch.
\end{enumerate}
\end{proof}

\begin{proof}
    \hypertarget{proof of cor: number of roots, leaves and branches}{(of Corollary \ref{cor: number of roots, leaves and branches})}
    \begin{enumerate}[(a)]
        \item Note that $|\{m\in\left(\Z/L\Z\right) \mid m=k\cdot\frac{L}{L_n}\}|=|\{k \mid 0\leq k< L_n\}|=L_n$. Therefore, by Theorem \ref{thm: zolotarev structure}, we have $L_n$ roots.

        \item Note that $|\{m\in\left(\Z/L\Z\right) \mid m=k\cdot(n,L)\}|=|\{k \mid 0\leq k< \frac{L}{(n,L)}\}|=\frac{L}{(n,L)}$. Then, by Theorem \ref{thm: zolotarev structure}, we have $\frac{L}{(n,L)}$ nodes which are not leaves, impliying that $|\text{roots}|=L-\frac{L}{(n,L)}$, and $|\text{branches}|=\frac{L}{(n,L)}-L_n$.

        \item 
    Let $F:=|\text{leaves}|$. By Lemma \ref{lem: (n,L) arrows}, the indegree of all nodes which are not leaves is $(n,L)$. Thus, there are exactly $(n,L)\cdot(L-F)$ edges in $Z(n,L)$. On the other hand, there are $L$ nodes, an therefore $L$ edges in total. Therefore, we have
    \[
    (n,L)\cdot(L-F)=L \iff F=L\cdot\frac{(n,L)-1}{(n,L)}
    \]
    \end{enumerate}

\end{proof}

\begin{proof}
    \hypertarget{proof of lem: m+kL_n is a root?} (of Lemma \ref{lem: m+kL_n is a root?})
By Theorem \ref{thm: zolotarev structure}, a node $r$ in $Z(n,L)$ is a root if, and only if, $\frac{L}{L_n}\mid r$. Then,
\[\frac{L}{L_n}\mid k\cdot L_n\iff k\equiv 0\Mod{\frac{L}{L_n}},\]
because $\left(\frac{L}{L_n},L_n\right)=1$.
Note that for each root $r$, we can construct a set of size $\frac{L}{L_n}$ given by $\{r+k\cdot L_n\mid0\leq k< \frac{L}{L_n}\}$. Moreover, the sets are disjoint, because, by construction, each set contains only one root.
Thus, we partitioned $Z(n,L)$ in $L_n$ sets.
\end{proof}

\subsection{Proofs of Theorem \ref{thm: number of roots with b_j}, Theorem \ref{thm: The Mother Tree} and Theorem \ref{thm: distance is invariant over +L_n}}

\begin{proof}
    \hypertarget{proof of thm: number of roots with b_j} (of Theorem \ref{thm: number of roots with b_j})
Note that $m$ is a root if, and only if, $m$ is in a cycle in $Z(n,L)$. Hence, there are  exactly $j\cdot b_j$ roots, where $b_j$ is number of cycles in $Z(n,L)$ of size $j$. It is clear that for any $j>L$, we have $b_j=0$. Therefore, by Observation \ref{Zolotarev graphs:formula for b_j}, the number of roots is given by:
\[
|\text{roots}|=\sum_{j=1}^{L}j\cdot b_j=\sum_{j=1}^{L}\sum_{d\mid j}\mu\left(\frac{j}{d}\right)\text{gcd}(n^d-1,L).
\]
\end{proof}

\begin{proof}
    \hypertarget{proof of thm: The Mother Tree}(of Theorem \ref{thm: The Mother Tree})
    It is easy to see that there exists an integer $k$ such that $n^km\equiv 0\mod L$ if, and only if, $L_n\mid m$. Also, $n^{d(L_n,0)}\cdot k\cdot L_n\equiv k\cdot0\equiv 0 \mod L$, which implies $d(L_n,0)\geq d(kL_n,0)$ for any integer $k$.
    
\end{proof}

\begin{proof}
    \hypertarget{proof of thm: distance is invariant over +L_n}(of Theorem \ref{thm: distance is invariant over +L_n})
Let $a$ and $b$ be two nodes, $i$ a root in the same connected component as $a$ and $b$ and suppose $a\equiv b \Mod{L_n}$. Let $b=a+kL_n$.
By Theorem \ref{thm: The Mother Tree}, for every node $kL_n\in\rm{Tree}(0)$, there exists an integer $\beta$ such that $n^{\beta}\cdot k L_n\equiv 0 \Mod{L}$. We have 
\[n^{d(b,i)}\cdot b=n^{d(a+kL_n,i)}\cdot (a+kL_n)\equiv i\Mod L.\]
Thus, there exists an integer $\alpha$ such that $d(b,i)+\alpha j>\beta$. Therefore, we obtain the following:

    \[
     n^{d(a+kL_n,i)+\alpha j}\cdot(a+kL_n)= n^{d(a+kL_n,i)+\alpha j}\cdot a+ n^{d(a+kL_n,i)+\alpha j}\cdot kL_n\equiv n^{d(a+kL_n,i)+\alpha j}\cdot a\equiv i \Mod L
    \]
    Then, for any integer $k$, the following holds:
    \[
    d(a+kL_n,i)+\alpha j\equiv d(a+kL_n,i)\equiv d(a,i) \Mod j
    \]

    Conversely, let $i$ be a root in a cycle of size $j$ and $a$ and $b$ two nodes in the same connected component as $i$ such that $d(a,i)\equiv d(b,i) \Mod j$.
    Let $\gamma:=b-a$ and $\alpha$ an integer such that $d(b,i)=d(a,i)+\alpha\cdot j$. Then, we have the following:
    \[
    i\equiv n^{d(b,i)}\cdot b=n^{d(b,i)}\cdot a+n^{d(b,i)}\cdot\gamma=n^{d(a,i)+\alpha j}\cdot a+n^{d(b,i)}\cdot\gamma \equiv i+n^{d(b,i)}\cdot\gamma \Mod L.
    \]
    Therefore, we conclude that $n^{d(b,i)}\cdot\gamma\equiv 0\Mod L$ and, by Theorem \ref{thm: The Mother Tree}, $\gamma$ is a multiple of $L_n$, hence $a\equiv b\Mod {L_n}$.
\end{proof}

\subsection{Proofs of Lemma \ref{lem: the max height} and Theorem \ref{thm: homogeneously height}}
\begin{proof}
\hypertarget{proof of lem: the max height} (of Lemma \ref{lem: the max height})
    Let $i$ be the root such that $n^{H(1)}\equiv i \Mod{L}$ and $k\in Z(n,L)$ be any node. We have
    \begin{align*}
        n^{H(1)}\cdot1&\equiv i\Mod{L}\\
        n^{H(1)}\cdot k&\equiv i\cdot k\Mod{L}\\
        \implies &H(k)\leq H(1),
    \end{align*}
    because $k\cdot i$ is a root in $Z(n,L)$.
\end{proof}

\begin{proof}
\hypertarget{proof of thm: homogeneously height} (of Theorem \ref{thm: homogeneously height})
Let $g:=${gcd}$(n,L)$ and suppose $Z(n,L)$ has homogeneous height $h$, i.e., all its leaves have the same height $h$. 
By Lemma \ref{lem: (n,L) arrows}, if $r$ is a root in $Z(n,L)$, its  indegree is $g$, from which $g-1$ edges connect $r$ to leaves and branches. Thus, there are $L_n\cdot(g-1)$ nodes $m$ with $H(m)=1$. Note that if $m$ is a branch, then its indegree is $g$. Therefore, we conclude, inductively, that, for every $i\leq h$, we have $|\{m\mid H(m)=i\}|=L_n\cdot(g-1)\cdot g^{i-1}$.

As there are $L$ nodes in total, we have the following:
\[
L=L_n+\sum_{m=1}^{h}L_n(g-1)g^{m-1}=L_n\left(1+(g-1)\sum_{m=0}^{h-1}g^{m}\right)
\]
\[
=L_n\left(1+(g-1)\cdot\frac{g^h-1}{g-1}\right)=L_n\cdot g^h,
\]
as desired.

Conversely, suppose $L_n=\frac{L}{g^h}$.
Let $L=\prod_{i=1}^{s} p_i^{\alpha_i}$, $n=\prod_{i=1}^{t} q_i^{\beta_i}$ and $\{r_i\mid 1\leq i\leq max\{s,t\}\}$ be the set of primes shared by $L$ and $n$. Then, $\alpha_i=c\cdot\beta_i$,  for some integer $c$. 

\noindent\textbf{Claim:} If $m$ is a leaf in $Z(n,L)$, then $H(m)=c$.

Indeed, by Theorem \ref{thm: zolotarev structure}, if $m$ is a leaf, then $g\nmid m$, i.e., gcd$(g,m)=\prod r_i^{\gamma_i}$. Moreover, there exists an integer $i'$ such that $\gamma_{i'}<\beta_{i'}$. But, by Theorem \ref{thm: zolotarev structure} $n^km$ is a root if, and only if, $\frac{L}{L_n}=g^h\mid n^km$, i.e., for all $i$, $\gamma_i+k\cdot \beta_i\geq h\cdot \beta_i$. But, because $\gamma_{i'}<\beta_{i'}$, we have $\gamma_{i'}+k\cdot \beta_{i'}\geq h\cdot \beta_{i'}\iff k\geq h$, which implies $H(m)=h$, concluding the proof.
\end{proof}

\subsection{Proof of Theorem \ref{thm: constructing isomorphism of subgraphs}, Corollary \ref{cor: mother tree isomorphic to Z(n,L)Z(L/L_n)}, Corollary \ref{cor: prune Z(n,L)} and Corollary \ref{cor: roots of Z(n,L)Z(n,L_n)}}

\begin{proof}
    \hypertarget{proof of thm: constructing isomorphism of subgraphs}(of Theorem \ref{thm: constructing isomorphism of subgraphs})
    The idea of the proof is the well-known identity:
\begin{equation}
\label{modulo identity}
a\equiv b\Mod c\iff\frac{a}{d}\equiv\frac{b}{d}\Mod{\frac{c}{d}},
\end{equation}
where  $d$ is any common divisor of $a, b, c$. 
Consider the function $\psi:G \rightarrow Z(n,\frac{L}{v})$, 
defined by:
\begin{align*}
\psi(x)=\frac{x}{v}.
\end{align*}
    Note that
    \begin{align*}
        \psi(x)\equiv \psi(y)\Mod{\frac{L}{v}}\iff \frac{x}{v}\equiv\frac{y}{v}\Mod{\frac{L}{v}}\iff x\equiv y\Mod{L}, 
    \end{align*} by \eqref{modulo identity}.
    Hence, $\psi$ is injective. Also, because the two graphs have $\frac{L}{v}$ nodes, we conclude that $\psi$ is a bijection. Moreover, $\psi$ has the isomorphism property \eqref{Isomorphic function}: If $a\longmapsto b\equiv na\mod L$, then by \eqref{modulo identity},
    \begin{align*}
        \psi(b)=\frac{b}{v}\equiv \frac{na}{v}=n\psi(a)\Mod {\frac{L}{v}}\iff \psi(a)\longmapsto \psi(b),
    \end{align*}
    concluding the proof.
\end{proof}

\begin{proof}
    \hypertarget{proof of cor: mother tree isomorphic to Z(n,L)Z(L/L_n)}(of Corollary \ref{cor: mother tree isomorphic to Z(n,L)Z(L/L_n)})
    By Theorem \ref{thm: The Mother Tree}, $x\in\tree(0)\iff L_n\mid x$. Thus, by Theorem \ref{thm: constructing isomorphism of subgraphs} the proof is done.
\end{proof}

\begin{proof}
    \hypertarget{proof of cor: prune Z(n,L)}(of Corollary \ref{cor: prune Z(n,L)})
    By Theorem \ref{thm: zolotarev structure}, all the non-leaves are multiples of $(n,L)$. Thus, it is a particular case of Theorem \ref{thm: constructing isomorphism of subgraphs}.
\end{proof}
\begin{proof}
    \hypertarget{proof of cor: roots of Z(n,L)Z(n,L_n)} (of Corollary \ref{cor: roots of Z(n,L)Z(n,L_n)})
    Apply Theorem \ref{thm: constructing isomorphism of subgraphs} with $v=\frac{L}{L_n}.$
\end{proof}

\subsection{Proof of Theorem \ref{thm: trees} and Corollary \ref{cor: number of nodes connected}}

\begin{proof}
    \hypertarget{proof of thm: trees}(of Theorem \ref{thm: trees})
    Here we use the notation of Algorithm Z, of section \ref{Algorithm:Zolotarev graphs}. 
     We proceed by backward induction, commencing with $L_n$, and ending with $L$.

    First, consider the Zolotarev Graph $Z(n,L_n)$. It is trivial that for each two roots $r_1\cdot\frac{L_n}{L}$ and $r_2\cdot\frac{L_n}{L}$ that we chose the theorem holds, because the subgraph $\tree(r)$ is just the node $r$.

    Now suppose the theorem holds for $Z(n,x_{j+1})$. Let show it hold for $Z(n,x_j)$. We know by Corollary \ref{cor: prune Z(n,L)} that $Z(n,x_{j+1})$ is isomorphic to the subgraph of the non-leaves of $Z(n,x_j)$. Now by Lemma \ref{lem: (n,L) arrows}, if a node in $Z(n,x_{j+1})$ has indegree $I$, then must connect exactly $(n,x_{j})-I$ new leaves in this node to obtain $Z(n,x_j)$. Thus it is clear that add this leaves preserve the isomorphism, therefore $\tree\left(r_1\cdot\frac{x_{j+1}}{L}\right)\simeq \tree\left(r_2\cdot\frac{x_{j+1}}{L}\right)\subset Z(n,x_{j+1})\implies\tree\left(r_1\cdot\frac{x_{j}}{L}\right)\simeq \tree\left(r_2\cdot\frac{x_{j}}{L}\right)\subset Z(n,x_{j})$, concluding the induction.
\end{proof}
\begin{proof}
    \hypertarget{proof of cor: number of nodes connected}(of Corollary \ref{cor: number of nodes connected}) By Theorem \ref{thm: trees}, we have that all the trees has the same number of nodes, and since by Theorem \ref{thm: The Mother Tree}, $\tree(0)$ has $\frac{L}{L_n}$ nodes, thus, if a connected component have a cycle of size $j$, then it must have exactly $j\cdot\frac{L}{L_n}$ nodes.
\end{proof}

%%%%%%%%%%%%%%%%%%%%%%%%%%%%%%%%%%%%%%%
%%%%%%%%%%%%%%%%%%%%%%%%%%%%%%%%%%%%%%%

\medskip
\section{Proofs for section \ref{sec:the eigenspaces}: The eigenspaces}

\subsection{Proofs of Theorem \ref{thm: basis for En(omega,L,k)}, Corollary \ref{cor: E_n(L)=E_n(L_n)} and Theorem \ref{thm:precise dimension of En(omega,L,k)}}
\begin{proof}
    \hypertarget{proof of thm: basis for En(omega,L,k)}(of Theorem \ref{thm: basis for En(omega,L,k)})  
First, we will show that
    $F_{\omega,L,\kappa,n,r}\in E_n(\omega,L,\kappa).$

 \noindent Indeed, if $i\neq r$ is connected to $r$, then $d(ni,r)=d(i,r)-1$. Thus, 
\begin{equation}
    a_{\omega,L,n,r}(ni)=\omega^{-d(ni,r)}=\omega^{-d(i,r)+1}=\omega a_{\omega,L,n,r}(i).
\end{equation}

If $i=r$, because the cycle has size $jm$ and $\omega$ is a $m$'th root of unity, we have
\begin{equation}
    a_{\omega,L,n,r}(nr)=\omega^{-d(ni,r)}=\omega^{-bm+1}=\omega=\omega a_{\omega,L,n,r}(r).
\end{equation}

\noindent If $i$ is not connected to $r$, then neither is $ni$, and therefore $a_{\omega,L,n,r}(ni)=\omega a_{\omega,L,n,r}(i)=0$.
\smallskip

\noindent If $r$ and $k$ share a cycle, $F_{\omega,L,\kappa,n,r}$ and $F_{\omega,L,\kappa,n,k}$ are linearly dependent, because $F_{\omega,L,\kappa,n,r}=\omega^{d(r,k)} F_{\omega,L,\kappa,n,k}$.
But taking just one root $r_i$ from each disjoint cycle $A_i\in B_m,1\leq i\leq |B_m|$ , we find that they are linearly independent, because 

\[
a_{\omega,L,n,r_j}(r_i)=
\begin{cases}
    1,&\text{ if $i=j$}\\
    0,&\text{ if $i\neq j$}
\end{cases}
\]

Then we can take in particular $\{
   F_{\omega, L, \kappa, n, r} \mid  \  \ r = \min(A) 
   \}$, implying that they are linearly independent

Now we show that they span the space. Let $f(x)=\sum_{k=0}^{\infty}k^{\kappa-1}a(k)x^k\in E_n(\omega,L,\kappa)$ and consider the Zolotarev graph $Z(n,L)$. By Corollary \ref{corollary:How to write an eigenfunction}, we have $a(k)=a(k+L)$, for all $k\in\Z_{\geq0}$ and $a(nk)=\omega a(k).$ If $k$ is a root in a cycle of size $j$, then $a(n^jk)=a(k)=\omega^ja(k)$ and, therefore, we have $a(k)=0$ or $\omega^j=1$. Indeed, if $\omega^j\neq 1$, then $a(k)=0$ and if $\omega^j=1$, we have $a(b)=a(k)\cdot a_{\omega,L,n,k}(b)$, for every node $b$ in the same connected component as $k$. Therefore,
\[f(x)=\sum_{\substack{{r=min(A)}\\
{A\in B_m}}}a(r)\cdot F_{\omega,L,\kappa,n,r}(x),\] concluding the proof.

\end{proof}
\begin{proof}
    \hypertarget{proof of cor: E_n(L)=E_n(L_n)}(of Corollary \ref{cor: E_n(L)=E_n(L_n)})
    Let $f\in E_n(\omega,L,\kappa)$. We will show that $f\in E_n(\omega,L_n,\kappa)$. Note that, by the definition of $E_n(\omega,L,\kappa)$, we have that  $level(f)\mid L$. However, by Lemma \ref{gcd(n, L)=1 for an eigenfunction}, gcd$(level(f),n)=1$. Therefore, $level(f)\mid L_n$ and we conclude that $E_n(\omega,L,\kappa)\subset E_n(\omega,L_n,\kappa)$.

    On the other hand, note that $L_n\mid L\implies E_n(\omega,L_n,\kappa)\subset E_n(\omega,L,\kappa)$, concluding the proof.
\end{proof}
\begin{proof}
    \hypertarget{proof of thm:precise dimension of En(omega,L,k)}(of Theorem \ref{thm:precise dimension of En(omega,L,k)})
    It follows from Theorem \ref{thm: basis for En(omega,L,k)} that $\dim E_n(\omega,L,\kappa)=|B_m|$, where $B_m$ is defined in \eqref{B_m}. Indeed, $|B_m|$ is the sum of the number of cycles of size $j$ in $Z(n,L)$, where $m\mid j$. Therefore, we have the following:
    \[|B_m|=
    \dim E_n(\omega,L,\kappa)=\sum_{m\mid j}^{L}b_j=\sum_{j=1}^{c}b_{mj}
    \]  
\end{proof}

\subsection{Proof of Theorem \ref{thm:precise dimension of Sn(L,k)}, Corollary \ref{cor: S_n(L)=R(L_n)} and Corollary \ref{cor: dim S_n multiplicative and periodic}}
\begin{proof}
    \hypertarget{proof of thm:precise dimension of Sn(L,k)}(of Theorem \ref{thm:precise dimension of Sn(L,k)})
    Let $A$ be a cycle of size $j$ in $Z(n,L)$, $r\in A$ and a $j$'th root of unity $\omega$. Accordingly, let $F_{\omega,L,\kappa,n,r}$ be a basic function of the non trivial eigenspace $E_n(\omega,L,\kappa)$. Note that the subspace spanned by $F_{\omega,L,\kappa,n,r}$ has dimension 1. Indeed, as there are $j$ $j$'th roots of unity, we have, for each cycle of size $j$, $j$ dimensions (one for each distinct root of unity). Thus,  $\dim S_n(L,\kappa)$ is equal the number of roots in $Z(n,L)$. Indeed, by Corollary \ref{cor: number of roots, leaves and branches}, we conclude that there are exactly $L_n$ roots and, therefore, $\dim S_n(L,\kappa)=L_n$, concluding the proof.
    
\end{proof}

\begin{proof}
    \hypertarget{proof of cor: S_n(L)=R(L_n)}(of Corollary \ref{cor: S_n(L)=R(L_n)})
    By Corollary \ref{cor: E_n(L)=E_n(L_n)}, we have that $E_n(\omega,L,\kappa)\subset\RR(L_n,\kappa)$, for every $\omega$. Thus, by definition, we conclude that $S_n(L,\kappa)\subset\RR(L_n,\kappa)$.
    Now, because of Theorem \ref{thm:precise dimension of Sn(L,k)} and Lemma \ref{lem: precise dimension of R_L,k},  $\dim S_n(L,\kappa)=\dim \RR(L_n,\kappa)$. Therefore, the two spaces are equal.
    
\end{proof}

\begin{proof}
    \hypertarget{proof of cor: dim S_n multiplicative and periodic}(of Corollary \ref{cor: dim S_n multiplicative and periodic})
    By Theorem \ref{thm:precise dimension of Sn(L,k)}, $\dim S_n(L,\kappa)=L_n$. Therefore, it is sufficient to show that $g(n,L):=L_n$ is periodic in $n$ and completely multiplicative in $L$.

    First, we show that $g$ is periodic in $n$. Note that
    \begin{align*}
        \gcd(n,L)=\gcd(n+L,L)\implies g(n,L)=g(n+L,L),
    \end{align*}
    by the algorithm described in example \eqref{recursive formula to L_n}.
    
    Now, we show that $g$ is completely multiplicative in $L$. Let $A$, $B$ and $n$ be integers. Consider the trivial decomposition $A=\frac{A}{A_n}\cdot A_n$ and $B=\frac{B}{B_n}\cdot B_n.$ We have
    \begin{align*}
        A\cdot B=\frac{A\cdot B}{A_n\cdot B_n}\cdot A_n\cdot B_n.
    \end{align*}
    Now, note that $(A_n\cdot B_n)\mid A\cdot B$ and it is indeed the largest divisor of $A\cdot B$ which is coprime to $n$, because if $p$ is a prime such that $p\mid \frac{A\cdot B}{A_n\cdot B_n}$, then $p\mid \frac{A}{A_n}$ or $p\mid \frac{B}{B_n}$. In any case, $p\mid n$, and we are done.
    \end{proof}
%%%%%%%%%%%%%%%%%%%%%%%%%%%%%%%%%%%%%%%
%%%%%%%%%%%%%%%%%%%%%%%%%%%%%%%%%%%%%%%
%%%%%%%%%%%%%%%%%%%%%%%%%%%%%%%%%%%%%%%
%%%%%%%%%%%%%%%%%%%%%%%%%%%%%%%%%%%%%%%

\medskip
\section{Proofs for section \ref{sec: Simultaneous eigenfunctions} Simultaneous eigenfunctions}

\medskip
\subsection{Proofs of Lemma \ref{lem: Defining C(f) and lambda_f}, and
Corollary \ref{cor: generating V_L}
}

\begin{proof} \hypertarget{proof of lem: Defining C(f) and lambda_f} (of Lemma \ref{lem: Defining C(f) and lambda_f})

    \begin{enumerate}[(a)]
        
        \item By Corollary \ref{corollary:How to write an eigenfunction}: \[
        U_n\left(f(x)\right)=\sum_{k=0}^{\infty}(nk)^{\kappa-1}a(nk)x^k=n^{\kappa-1}\omega_f(n)\cdot f =\sum_{k=0}^{\infty}(nk)^{\kappa-1}\omega_f(n)a(k)x^k \Longleftrightarrow a(nk)=\omega_f(n)a(k),
        \] for all $k\geq 0$.
        
        \item Let $f\in\RR$. Then, 
        \[
        U_i\circ U_j (f)= U_i(\sum_{k=0}^{\infty} a(j\cdot n) x^k)=\sum_{k=0}^{\infty} a(i\cdot j\cdot k) x^k=U_{i\cdot j}.
        \]
        Thus, $U_i\circ U_j=U_{i\cdot j}$.
        Therefore, $U_n\circ U_m(f)=U_n(\omega_f(m)f)=\omega_f(n)\omega_f(m)f=U_{mn}(f)$, implying $mn\in\mathcal{C}(f)$ and $\omega_f(m\cdot n)=\omega_f(m)\cdot\omega_f(n)$.
        
        \item Let $f\in E_m(\omega_f(m),L,\kappa)$. By Lemma \ref{corollary:How to write an eigenfunction}, $a(k+L)=a(k)$ for all $k\geq 0$. Applying the operator $U_{m+L}$, we have:
        \begin{align*}
        U_{m+L}\left(f(x)\right)
        &=\sum_{k=0}^{\infty}((m+L)k)^{\kappa-1}a((m+L)k)x^k \\
        &=\left(m+L\right)^{\kappa-1}\sum_{k=0}^{\infty}(k)^{\kappa-1}a(mk)x^k\\
        &=\left(m+L\right)^{\kappa-1}\omega_f(m)\sum_{k=0}^{\infty}(k)^{\kappa-1}a(k)x^k,
        \end{align*}
        which implies $m+L\in\mathcal{C}(f)$ and $\omega_f(m+L)=\omega_f(m)$.
        
        \item By the latter two proofs, we easily conclude that $\omega_f$ is a completely multiplicative and periodic function, with period $L$. Moreover, by Lemma \ref{gcd(n, L)=1 for an eigenfunction}, if $\left(d,L\right)>1$ and $d\in\mathcal{C}(f)$, implying $\omega_f(d)=0$.
        Therefore, there exists a Dirichlet character $\chi$ modulo $L$  such that $\omega_f=\chi|_{\mathcal{C}(f)}$, concluding the proof of the Lemma.
    \end{enumerate}
\end{proof}

\begin{proof} \hypertarget{proof of cor: generating V_L}{} (of Corollary \ref{cor: generating V_L})
    If Level$(f)=\overline{L}$, then, by Lemma \ref{lem: Defining C(f) and lambda_f}, there exists a Dirichlet character modulo $\overline{L}$ such that $\chi(n)=\omega_f(n)$, for all $n\in\mathcal{C}(f)$. In particular, for $\mathcal{C}(f)=\Z$, the result follows.
\end{proof}

\begin{proof}
\hypertarget{proof of (M+kc,L)=1} (of Lemma \ref{(M+kc,L)=1})
Suppose $p\mid c$ and that there exists $k\in\Z$ with $\left(M+k\cdot c,L\right)>1$. Then, there exists a prime $q$ such that $q \mid L$ and $q \mid M+k\cdot c$. As $L=p\cdot c$, we have $q\mid p$ or $q\mid c$. But $q\mid p \mid c$, then $q\mid c$. Therefore:
\[
\begin{cases}
    q\mid c \\
    q\mid M+k\cdot c
\end{cases}\Rightarrow q\mid M,
\]
contradicting $(M,L)=1$.

Conversely, assume that for all $k\in\Z$, $(M+kc,L)=1$.
Suppose that  $p\nmid c$. Then, there exists is a prime $q$ such that  $q\mid p$ and $q\nmid c$. But $p\mid L \implies q\mid L$, and $(M,L)=1 \implies q\nmid M$.

Now, if $(M+kc,L)=1$, then $q\nmid M+kc$.
Let $\overline{M}$ and $\overline{q}$ be such that $M\equiv \overline{M} \Mod q$ and $c\equiv \overline{c}\Mod q$. However, $\left(\Z/q\Z\right)^*$ is a field, so there exists $(\overline{c})^{-1}$, the inverse of $\overline{c}$. Let $k^*=(\overline{c})^{-1}\cdot (-\overline{M})$. We have $\overline{M}+\left((\overline{c})^{-1}\cdot (-\overline{M})\right)\overline{c}\equiv 0 \pmod q$, which means $k^*$ is such that $q\mid M+k^*c$, contradicting $(M+kc,L)=1$ for all $k$.

\end{proof}

\subsection{Proofs of Lemma \ref{lem: Defining A_l} and Theorem \ref{thm:precise dimension of V_L,k}}

\begin{proof}
    \hypertarget{proof of lem: Defining A_l} (of Lemma \ref{lem: Defining A_l})
    \begin{enumerate}[(a)]
        \item Suppose $a\in A_L$. Then, we must have $a=\prod_{j\in S}q_j$ and, therefore,$\frac{L}{a}=\prod_{j\in S^c}q_j$. 
        Thus, $a\mid L$ and $\left(a,\frac{L}{a}\right)=1$.

        Conversely, suppose $a\mid L$ and $\left(a,\frac{L}{a}\right)=1.$ Then, 
        \[a=\prod_{j=1}^{n}p_j^{\beta_j}\Rightarrow \frac{L}{a}=\prod_{j=1}^{n}p_j^{\alpha_j-\beta_j}.\] 
        Thus, 
        \[1=\left(a,\frac{L}{a}\right)=\prod_{j=1}^{n}p_j^{min\{\alpha_j-\beta_j,\beta_j\}}\iff \forall j\text{, }\beta_j=0 \text{ or } \beta_j=\alpha_j.\]
        In any case, we have $a\in A_L$, concluding the proof of (a).

        \item Let $C=\{c\; \mid \; L\equiv0\pmod{c} \text{ and } \forall a>1\text{ such that }c\equiv 0 \pmod a\Rightarrow a\cdot c\nmid L\}$. First, we show that $A_L\subset C$. Let $k\in A_L$. Then, $k=\prod_{j\in S}q_j$. Note that $p_i\mid k$ if, and only if $i\in S$, which implies $p_i\cdot k\nmid L$.
        Conversely, we show that $C\subset A_L$. Let $c\in C$ and $c=\prod_{j\in S}p_j^{\beta_j}$.
         If there exists an integer $j$ such that $\beta_j<\alpha_j$, where $\alpha_j$ is such that $q_j=p_j^{\alpha_j}$, we can choose $a=p_j$, implying $c\cdot a\mid L$, a contradiction, concluding the proof of (b).
        \item Let $\mathcal{J}_n=\{1,\dots, n\}$. Note that $\mathcal{P}(\mathcal{J}_n)$
    can be partitioned in the following manner:
    \[
\mathcal{P}(\mathcal{J}_n)=\left\{S\mid S\subset \{1, \dots n\}\setminus\{i\}\right\}\cup\left\{S\cup \{i\}\mid S\subset \{1, \dots n\}\right\}
    \]
    Let $\mathcal{A}=\left\{S\mid S\subset \{1, \dots n\}\setminus\{i\}\right\}$ and $\mathcal{B}=\left\{S\cup \{i\}\mid S\subset \{1, \dots n\}\right\}$. 
    
    If $S\in\mathcal{A}$, then $a=\prod_{j\in S}q_j\in A_{\frac{L}{q_i}}$. On the other hand, if $S\in\mathcal{B}$, then $b=\prod_{j\in S}q_j\in q_i\cdot A_{\frac{L}{q_i}}$, concluding the proof of (c).
        \item Note that $|A_L|=|\mathcal{P}(\mathcal{J}_n)|$, and we are done.
     
        \item Let $\mathcal{L}=\{L_n\mid\forall n\in\N\}$. First, we show that $A_L\subset \mathcal{L}$. Let $a\in A_L$. By item (a), we have that $a$ and $\frac{L}{a}$ are coprime. Indeed, $a$ is the largest divisor of $L$ coprime with $\frac{L}{a}$.  Thus, $a=L_{\frac{L}{a}}\Rightarrow A_L\subset \mathcal{L}$.
        Now, we show that $\mathcal{L}\subset A_L$. Suppose $l\in\mathcal{L}$. Then there exists an integer $m$ such that $l=L_m$.
        Now, because $l\mid L$, we can write $l=\prod_{j=1}^{n} p_j^{\beta_j}$, where $0\leq\beta_j\leq\alpha_j$. Indeed, it is clear that if $(p_i,m)>1$, then $\beta_i=0$, and if $(p_i,m)=1$, then $\beta_i=\alpha_i$. Indeed, if not, then $p_i$ is not the greatest divisor of $L$ coprime with $m$. Thus we conclude $l=\prod_{j\in S}q_j\in A_L$, and we are done.
    \end{enumerate}
\end{proof}

\begin{proof}
    \hypertarget{proof of thm:precise dimension of V_L,k} (of Theorem \ref{thm:precise dimension of V_L,k})

\noindent
{\bf Step 1}. \  We first show that 
$\bigcup_{M\in A_L}\mathfrak{B}_{M,\kappa}$
    spans $V(L,\kappa)$.
    By Corollary \ref{cor: generating V_L} we have for $f\in V(L,\kappa)$,
    $a(n\cdot k)=\chi(n)a(k)$, in particular $a(n)=\chi(n)a(1)$. Letting $c=a(1)$, we have $a(n)=c\cdot\chi(n)$, which proves that $\bigcup_{M\mid L}\mathfrak{B}_{M,\kappa}$ spans all the space.
    
    Next, if we have $L=p\cdot c$, where $p\mid c$, and $p$ is a prime. Now, because of the orthogonality of the Dirichlet characters, we have
    \[
    \left<\chi \text{ mod c}\right>=\left<1_{a \text{ mod c }},0<a\leq c\text{ and }(a,c)=1\right>.
    \]
    
    And since $p\mid c$, then $(a,L)=1\iff(a,c)=1$. Then by Lemma \ref{(M+kc,L)=1} we have that for all $0<a\leq c$ coprime with $c$, $(a+k\cdot c,L)=1$, and also:
    \[
    1_{a \text{ mod c }}=\sum_{k=0}^{p-1}1_{a+k\cdot c \text{ mod L}},
    \]
   which proves that
    \[
    \left<1_{a \text{ mod }c},0<a\leq c\text{ and }(a,c)=1\right>\subset \left<1_{a \text{ mod }L},0<a\leq L\text{ and }(a,L)=1\right>.
    \]
    In other words,
    \[
    \left<\mathfrak{B}_{c,\kappa}\right>\subset\left<\mathfrak{B}_{L,\kappa}\right>
    \]
    Now, by Lemma \ref{lem: Defining A_l} (b) we have 
    \[
    \left<\bigcup_{M\mid L}\mathfrak{B}_{M,\kappa}\right>=\left<\bigcup_{M\in A_L}\mathfrak{B}_{M,\kappa}\right>=V\left(L,\kappa\right)
    \]

\noindent\textbf{Step 2}.   Now, we will proceed by induction in the number of prime divisors of $L$ to show that this set is also linear independent. For the base case, let $L=p^v,$ where $p$ is a prime.
    We know, by the start of the proof, that $\mathfrak{B}_{p^v,\kappa}\cup\mathfrak{B}_{1,\kappa}$ generates $V(p^v,\kappa)$, and we now show that it is in fact a basis for $V(p^v,\kappa).$

    Consider the collection of functions $\left\{f_i\mid i=1,\dots, \phi(p^v)\right\}=\mathfrak{B}_{p^v,\kappa}$ defined by:
    \[
    f_i(x)=\sum_{k\geq0}\chi_i(k)k^{\kappa}x^k,
    \]
where $\chi_0(k)=1$, for all integers $k$.
The other $\chi_i$'s  are characters modulo $p^v$. 
Suppose we have any vanishing linear combination of the $f_i$'s:
\[
0=\sum_{i=0}^{\phi(p^v)}\alpha_i f_i(x)
=\sum_{i=0}^{\phi(p^v)}\sum_{k\geq0}\alpha_i\chi_i(k)k^{\kappa}x^k
=\sum_{k\geq0}\left(\sum_{i=0}^{\phi(p^v)}\alpha_i\chi_i(k)\right)k^{\kappa}x^k
 \iff \sum_{i=0}^{\phi(p^v)}\alpha_i\chi_i(k)=0,
\]
    for every $k\geq0$. In particular, if $k=p^v$, then
\[
\sum_{i=0}^{\phi(p^v)}\alpha_i\chi_i(p^v)
=\alpha_0\chi_0(p^v)=\alpha_0
=0.
\]
Therefore, for all $k\geq 0$,
\[
\sum_{i=0}^{\phi(p^v)}\alpha_i\chi_i(k)
=\sum_{i=1}^{\phi(p^v)}\alpha_i\chi_i(k)
=0\iff \alpha_i=0 \text{ for all } i=1,\dots,\phi(p^v),
\]
    using the orthogonality of the Dirichlet characters. Thus, we have 
    \[
    \dim V(p^v,\kappa)=\phi(p^v)+1.
    \]

    \noindent{\bf Step 3}. Now, for the induction step, we will show that
    \begin{equation*}
\bigcup_{M\in A_L}\mathfrak{B}_{M,\kappa} \text{ is a basis for } V(L,\kappa).
\end{equation*} 
Assume it holds for $c=p_1^{v_1}p_2^{v_2}\cdots p_N^{v_N}$, the prime factorization of $c$, and for $L=c\cdot p^v$, where $p\neq p_i,$ for all $i\leq N$. Let $\chi_{m}(k;i)$  denote  the $i$'th character modulo $m$. 
    Then the set that span is a basis for $V(L,\kappa)$ if, and only if, all $f_{i,M}$ are linearly independent, where

    \[
    f_{i,M}(x)=\sum_{k\geq0}\chi_{M}(k;i)k^{\kappa}x^k.
    \]
Suppose we have any vanishing linear combination of the $f_{i,M}$´s:
    
    \begin{align*}
    0
    &=\sum_{M\in A_L}\sum_{i=1}^{\phi(M)}\beta_{i,M}\cdot f_{i,M}(x)\\
    &=\sum_{M\in A_L}\sum_{i=1}^{\phi(M)}\beta_{i,M}\cdot\sum_{k\geq0}\chi_{M}(k;i)k^{\kappa}x^k\\
    &=\sum_{k\geq0}\left(\sum_{M\in A_L}\sum_{i=1}^{\phi(M)}\beta_{i,M}\cdot\chi_{M}(k;i)\right)k^{\kappa}x^k\\
    &\iff\sum_{M\in A_L}\sum_{i=1}^{\phi(M)}\beta_{i,M}\cdot\chi_{M}(k;i)=0 \text{ for all } M\in A_L, i=1,\dots,\phi(M)\text{ and }\beta_{i,M}\in\C.
    \end{align*}

    By Lemma \ref{lem: Defining A_l}, we know that $A_L=(p^v+1)A_c$, and we note that $A_c\cap p^vA_c=\varnothing.$ Therefore, if $f_{i,M}(x)=\sum_{k\geq0}a(k)k^{\kappa}x^k$, we have
    \begin{align*}
    a(k)=&\sum_{M\in A_L}\sum_{i=1}^{\phi(M)}\beta_{i,M}\cdot\chi_{M}(k;i)\\
    =&\sum_{M\in A_c}\sum_{i=1}^{\phi(M)}\beta_{i,M}\cdot\chi_{M}(k;i)+\sum_{M\in p^vA_c}\sum_{i=1}^{\phi(M)}\beta_{i,M}\cdot\chi_{M}(k;i)\\
    =&\sum_{M\in A_c}\sum_{i=1}^{\phi(M)}\beta_{i,M}\cdot\chi_{M}(k;i)+\sum_{M\in A_c}\sum_{i=1}^{\phi(M\cdot p^v)}\beta_{i,M\cdot p^v}\cdot\chi_{M\cdot p^v}(k;i)
    \end{align*}

        Now, if $k=n\cdot p^v$, $n\in \{1, 2, \dots, c$\}, then it is clear that 
    \[
    a(k)=\sum_{M\in A_c}\sum_{i=1}^{\phi(M)}\beta_{i,M}\cdot\chi_{M}(k;i).
    \]
    But because $(c,p)=1$, then $k=n\cdot p^v\mod c$ is a permutation of the set $\{1,2,\dots,c\}$, and by hypothesis, for all $k\in\Z/c\Z$,
    \[
    a(k)=\sum_{M\in A_c}\sum_{i=1}^{\phi(M)}\beta_{i,M}\cdot\chi_{M}(k;i)=0\iff\beta_{i,M}=0 \text{ for all } M\in A_c \text{ and } i=1,\dots,\phi(M)
    \]
We know that for any $k\in \Z/c\Z$, we have:
    \[
    a(k)=\sum_{M\in A_c}\sum_{i=1}^{\phi(M\cdot p^v)}\beta_{i,M\cdot p^v} \cdot \chi_{M\cdot p^v}(k;i).
    \]
    
    If $h=M\cdot p^v$, then for all $j\in\{1,\dots,\phi(p^v)\}$ and $l\in\{1,\dots,\phi(M)\}$, there is an unique correspondent $i\in\{1,\dots,\phi(h)\}$ such that $\chi_{p^v}(k;j)\cdot\chi_{M}(k;l)=\chi_{h}(k;i)$, using the isomorphism of a finite abelian group and its character group. Then:
    
    \begin{align} 
    a(k)=&\sum_{M\in A_c}\sum_{i=1}^{\phi(M\cdot p^v)}\beta_{i,M\cdot p^v}\cdot\chi_{M\cdot p^v}(k;i)\nonumber\\
    =&\sum_{M\in A_c}\sum_{i=1}^{\phi(M\cdot p^v)}\beta_{i,M\cdot p^v} \cdot \chi_{M}(k;l_i)\cdot\chi_{p^v}(k;j_i).\label{character sum}
    \end{align}
    
    Now, reorganizing \eqref{character sum} and letting $\gamma_{l_i,j_i,M}:=\beta_{i,M\cdot p^v}$, using the same uniqueness relation between $l_i,j_i$ and $i$, we have:
    \[
    \sum_{M\in A_c}\sum_{l=1}^{\phi(M)}\chi_{M}(k;l)\sum_{j=1}^{\phi(p^v)}\gamma_{l,j,M}\cdot\chi_{p^v}(k;j).
    \]
    
Using the orthogonality of Dirichlet characters, we observe that
\[
    \left<\chi\text{ mod }p^v\right>=\left<1_{a\text{ mod }p^v}\text{, }0<a\leq p^v\text{ and }(a,p^v)=1\right>,
\]
and hence we have
$\sum_{j=1}^{\phi(p^v)}\gamma_{l,j,M}\cdot \chi_{p^v}(k;j)=\sum_{j=1}^{\phi(p^v)}\alpha_{l,j,M}\cdot 1_{l_j\text{(mod $p^v$)}}(k)$, where $(l_j,p^v)=1$. Then:
    \begin{align*}
        &\sum_{M\in A_c}\sum_{l=1}^{\phi(M)}\chi_{M}(k;l)\sum_{j=1}^{\phi(p^v)}\gamma_{l,j,M}\cdot\chi_{p^v}(k;j)\\
        =&\sum_{M\in A_c}\sum_{l=1}^{\phi(M)}\chi_{M}(k;l)\sum_{j=1}^{\phi(p^v)}\alpha_{l,j,M}\cdot 1_{l_j\text{(mod $p^v$)}}(k)
    \end{align*}

    Now, let $k=l_t+i\cdot p^v$, $i=0,\dots,c-1$, with $t=1,\dots,\phi(p^v)$. Then:
    \[
    a(k)=\sum_{M\in A_c}\sum_{l=1}^{\phi(M)}\chi_{M}(k;l)\cdot\alpha_{l,t,M}
    \]
    Note that $i\cdot p^v$ is a permutation in $\Z/c\Z$ and, therefore, $l_t+i\cdot p^v$ is also a permutation. Thus, for all $t=1,\dots \phi(p^v)$ and $k\in\Z/c\Z$:
    \[
    a(k)=\sum_{M\in A_c}\sum_{l=1}^{\phi(M)}\alpha_{l,t,M}\cdot\chi_{M}(k;l)=0\iff\alpha_{l,t,M}=0 \text{ for all } M\in A_c \text{ and } l=1,\dots,\phi(M),
    \]
    and we are done.
\end{proof}

%%%%%%%%%%%%%%%%%%%%%%%%%%%%%%%%%%%%%%%
%%%%%%%%%%%%%%%%%%%%%%%%%%%%%%%%%%%%%%%
%%%%%%%%%%%%%%%%%%%%%%%%%%%%%%%%%%%%%%%
%%%%%%%%%%%%%%%%%%%%%%%%%%%%%%%%%%%%%%%

\medskip
\section{Proofs for section \ref{sec: The Kernel}: The Kernel}
\subsection{Proofs of Lemma \ref{lem: equivalence for the kernel}, Corollary \ref{cor:precise dimension of ker(U_n)} and Theorem \ref{thm: diagonalizable of U_n} }

\begin{proof} \hypertarget{proof of lem: equivalence for the kernel} (of Lemma \ref{lem: equivalence for the kernel})    
Suppose $f\in\ker(U_n)$, i.e., $a(nk)=0,$ for all $1\leq k\leq L.$ This is precisely equivalent to $a(b)=0$ for every node $b$ which is not a leaf, because, by definition,  if $b$ is a leaf, then there exists an integer $k$ such that $nk\equiv b\pmod L$.

Conversely, we will show that if $a(k)=0$ for every non-leaf $k$, then $f\in\ker(U_n)$.
    Let 
    \[
    f(x)=\sum_{k=0}^{\infty}k^{\kappa-1}a(k)x^k.
    \] We have
    \[
    U_n(f(x))=\sum_{k=0}^{\infty}(nk)^{\kappa-1}a(nk)x^k.
    \]
    By definition, $nk$ is not a leaf in $Z(n,L)$, thus $a(nk)=0$ for all $k$ and therefore $f\in\ker(U_n).$ 

\end{proof}
\begin{proof}
    \hypertarget{proof of cor:precise dimension of ker(U_n)}(of Corollary \ref{cor:precise dimension of ker(U_n)})
    By Corollary \ref{cor: number of roots, leaves and branches}, there are $\frac{L}{(n,L)}$ leaves in $Z(n,L)$, and we are done.
\end{proof}
\begin{proof}
    \hypertarget{proof of thm: diagonalizable of U_n} (of Theorem \ref{thm: diagonalizable of U_n})
    \begin{enumerate}
    \item ($a\iff b$) $U_n$ is diagonalizable $\iff$ $\dim S_n+\dim \ker U_n=\dim\RR(L,\kappa)$. We have by Lemma \ref{lem: precise dimension of R_L,k}, Theorem \ref{thm:precise dimension of Sn(L,k)} and Corollary \ref{cor:precise dimension of ker(U_n)} that 
    \[L_n+L-\frac{L}{(n,L)}=L\implies L_n=\frac{L}{(n,L)}.\]
    
    \item ($b \iff c$) Follows imediatly from Theorem \ref{thm: zolotarev structure}.
    
    \item ($c \implies d$) It is clear that $1$ is a leaf, and $1\longmapsto n$, meaning that $n$ is a root.
    
    \item ($d \implies e$) It is clear that $1\longmapsto n$, then $H(1)=1$.
    
    \item ($e \implies c$) By Lemma \ref{lem: the max height}, the height of the Zolotarev graph is $h=H(1)=1$ and, therefore, the leaves are directly connected to the roots.
    \end{enumerate}

\end{proof}

\bigskip
\section{Further Remarks and questions}
\begin{rem}
Suppose we fix a denominator $B(x)$ of a potential eigenfunction $f$ of $U_n$, with a fixed eigenvalue $\lambda$.  What is the collection of all numerator polynomials $A(x)$ such that $\frac{A(x)}{B(x)}$ is an eigenfunction of $U_n$?
\end{rem}

In general there is not a unique eigenfunction $\frac{A(x)}{B(x)}$, associated with the data above. An example of this phenomena is given by   $f_{6,4}$ and $f_{6,5}$ in Appendix $9$ of \cite{GilRobins})
(see also Corollary 3.12 of \cite{GilRobins}).

\begin{question}
What is the Kernel of $U_n:\RR\rightarrow\RR$?
\end{question}

\begin{question}
Can we make further progress on the Artin conjecture by analyzing carefully the eigenspaces $E_n(\omega,p,\kappa)$?
\end{question}
%%%%%%%%%%%%%%%%%%%%%%%%%%%%%%%%%%%%%%%%%%%%
%%%%%%%%%%%%%%%%%%%%%%%%%%%%%%%%%%%%%%%%%%%%

\appendix

\section{Alternative proofs}

\subsection{Alternative proofs of Corollary \ref{cor: roots of Z(n,L)Z(n,L_n)}, Theorem \ref{thm:precise dimension of En(omega,L,k)} and Corollary \ref{cor: E_n(L)=E_n(L_n)}}

\begin{proof}
\hypertarget{annex proof of cor: roots of Z(n,L)Z(n,L_n)}(of Corollary \ref{cor: roots of Z(n,L)Z(n,L_n)})
    First, we show that for any integer $d$, $(n^d-1,L)=(n^d-1,L_n)$.
As $(n^d-1,n^d)=1$, $(n^d-1,L)$ does not have any common factor with $n^d$.
The same argument shows that $(n^d-1,\frac{L}{L_n})=1$. Thus,

\[(n^d-1,L)=(n^d-1,L_n\frac{L}{L_n})=(n^d-1,L_n)\]

    The second step is to show that two subgraphs $G$ and $H$ of $Z(n,L)$ consisting only of roots are isomorphic if, and only if, all the $b_j$´s associated to that cycle are identical.

    Indeed, if there exists a function that preserves adjacency between nodes, then it trivially preserves cycle sizes. The converse is also clear.

    Now, note that the $b_j$ associated with $Z(n,L)$ are given by 

    \[
    b_j =  \frac{1}{j}\sum_{d\mid j} \mu\left(\frac{j}{d}\right)\gcd(n^d-1, L)=\frac{1}{j}\sum_{d\mid j} \mu\left(\frac{j}{d}\right)\gcd(n^d-1, L_n)
    \]
    Indeed this is exactly the formula for $b_j$ associated to $Z(n,L_n)$.
\end{proof}

\begin{proof}
\hypertarget{anex proof of thm:precise dimension of En(omega,L,k)}
(of Theorem \ref{thm:precise dimension of En(omega,L,k)})
    Let $(n,L)=1$. By Theorem \ref{thm:eigenvalues and kappa} and Lemma \ref{gcd(n, L)=1 for an eigenfunction}, if $f\in E_n(\omega,L,\kappa)$, then $U_n(f)=n^{\kappa-1}\omega f$, where $\omega$ is a root of unity.
    Indeed, there exists a unique integer $m$ such that $\omega$ is a primitive $m$'th root of unity.

    Now, let $A$ be a cycle of size $j$ in $\tau(n,L)$. If $k\in A$, then $a(k)=a(kn^j)=\omega a(kn^{j-1})=\omega^ja(k)$, which implies $\omega^j=1$ or $a(k)=0$. But, if  $m\mid j$, then $\omega^j=1$, because $\omega$ is a primitive $m$´th root of unity. Hence, we have that $a(k)$ is a free variable. Indeed, once $a(k)$ is fixed, $a(kn^i)=\omega^ia(k)$, for all $b\in A$, with $b=k\cdot n^i$, for some integer $i$.
    Now, if $m\nmid j$, then $\omega^j\not=1$ and, therefore, $a(k)=0$, implying $a(k n^i)=0$, for all $i$.

    Thus, we conclude that $\dim E_n(\omega,L,\kappa)=\sum_{m\mid j}^{c} b_j$, where $c:=\ord_L(n)$.
    Now, let $L$ be any integer, not necessarily coprime to $n$. Therefore, by Corollary \ref{cor: E_n(L)=E_n(L_n)}, $E_n(\omega,L,\kappa)=E_n(\omega,L_n,\kappa)$, concluding the proof.

\end{proof}

\subsection{Alternative proof of Theorem \ref{thm:precise dimension of Sn(L,k)}}
\begin{proof}
\hypertarget{anex proof of thm:precise dimension of Sn(L,k)}
    (of Theorem \ref{thm:precise dimension of Sn(L,k)})

Note that
    \begin{equation}
    \label{eq: dimension S_n}
        \dim S_n(L,\kappa)=\sum_{\omega\in\C \setminus \{0\}}\dim E_n(\omega,L,\kappa)
    \end{equation}

    Now, let $\left(n,L\right)=1$. Then, by Theorem \ref{thm:eigenvalues and kappa} and Lemma \ref{gcd(n, L)=1 for an eigenfunction}, $\omega$ is a $c$'th root of unity, where $c=\ord_L(n)$.

    Moreover, $\dim E_n(\omega,L,\kappa)$ only depends on $n$, $L$ and $M$, where $\omega$ is a primitive $M$'th root of unity (and, therefore, $M\mid c$. Note that the fact that there exist $\phi(M)$ primitive $M$'th roots of unity implies that we can rewrite \eqref{eq: dimension S_n} as
    
    \[
    \dim S_n(L,\kappa)=\sum_{M\mid c}\phi(M)\dim E_n(e^{\frac{2\pi i}{M}},L,\kappa)
    \]
    Now, applying Theorem \ref{thm:precise dimension of En(omega,L,k)}, we have    \[
    \dim S_n(L,\kappa)=\sum_{M\mid c}\phi(M)\sum_{M\mid j}^{c} b_j
    =\sum_{M\mid c}\sum_{M\mid j}^{c} \phi(M) b_j.
    \]
    Now, by Lemma \ref{lem:defining b_j}, if $j\nmid c$, then $b_j=0$. I.e.,
    \[
    \sum_{M\mid c}\sum_{M\mid j}^{c} \phi(M) b_j=\sum_{M\mid c}\sum_{M\mid j\mid c}\phi(M) b_j.
    \]
    Alternatively, if $c=\prod_{i=1}^{k}p_i^{\alpha_i}$, then
    $m=\prod_{i=1}^{k}p_i^{\beta_i}$ and
    $j=\prod_{i=1}^{k}p_i^{\gamma_i}$,
    where $\beta_i\leq\gamma_i\leq\alpha_i$. Now, it is easy to see that we can rewrite the sum as
    \[
    \sum_{M\mid c}\sum_{M\mid j\mid c}\phi(M) b_j=\sum_{j\mid c}\sum_{M\mid j} \phi(M) b_j
    =\sum_{j\mid c} b_j\sum_{M\mid j} \phi(M)
    =\sum_{j\mid c} j b_j.
    \]
    Now, because $j\nmid c\Rightarrow b_j=0$, we have
    \[
    \sum_{j\mid c} j b_j=\sum_{j\mid c} j b_j + \sum_{\substack{{j\nmid c}\\{1\leq j<c}}}j b_j =\sum_{j=1}^{c} j b_j.
    \]
     By Lemma \ref{lem:defining b_j}, it follows that
    \[
    \sum_{j=1}^{c} j b_j=L,
    \]
    showing that if $(n,L)=1$, then $\dim S_n(L,\kappa)=L$.
    
    Now, let $n$ and $L$ be any positive integers. By Corollary \ref{cor: E_n(L)=E_n(L_n)}, $E_n(\omega,L,\kappa)=E_n(\omega,L_n,\kappa)$ and,  therefore, $S_n(L,\kappa)=S_n(L_n,\kappa)$, concluding the proof.
\end{proof}

\section{Zolotarev graph examples, with L=20}
\label{appendix: examples}
Here we exhibit all  Zolotarev graphs $Z(n, L)$, with $L=20$, and $2\leq n \leq 18$, with $\gcd(n, 20)>1$.

\begin{figure}[H]
\includegraphics[width=8cm]{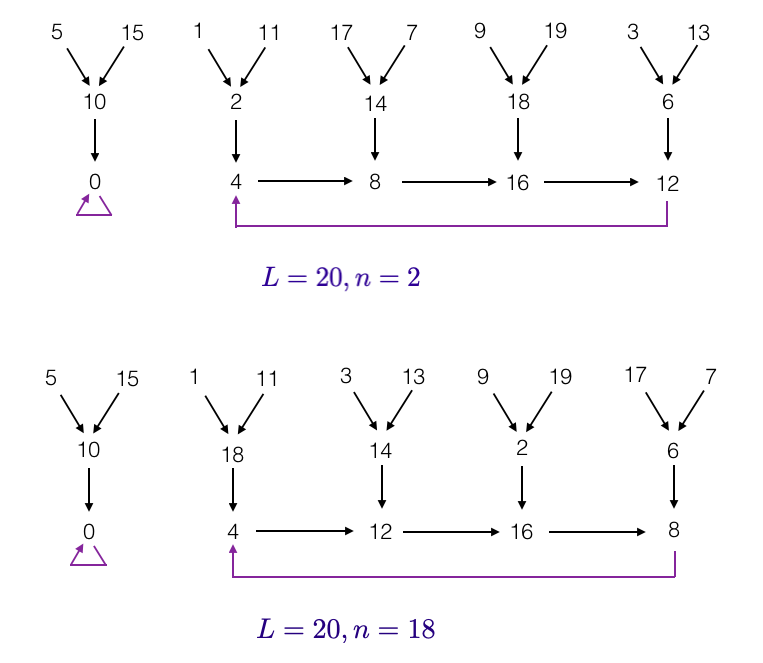}
\caption{ Each of these graphs has $10$ leaves, $5$ branches, $5$ roots, and $2$ connected components.} 
\label{newL=20,n=2,18}
\end{figure}

\begin{figure}[H]
\includegraphics[width=9cm]{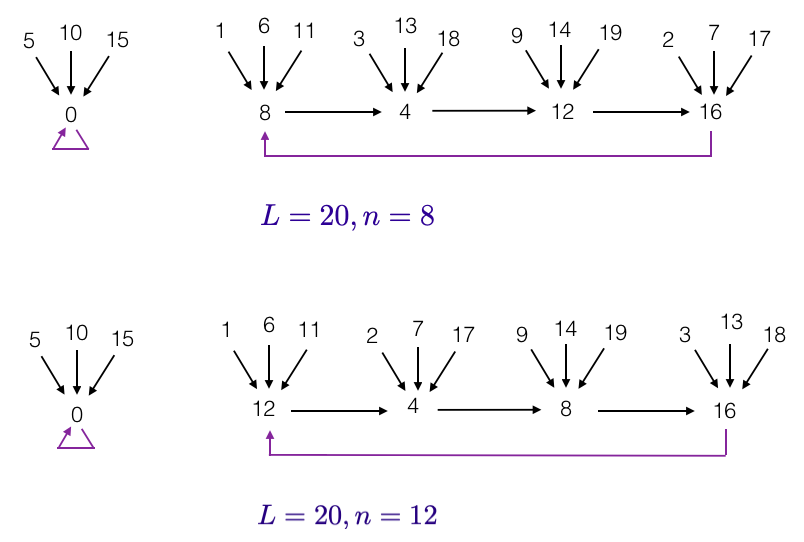}
\caption{The two Zoloterev graphs $Z(8, 20)$ and $Z(12, 20)$ are isomorphic. Each of them has $15$ leaves, no branches, $5$ roots, and $2$ connected components. } 
\label{newL=20,n=8,12}
\end{figure} 

\begin{figure}[H]
\includegraphics[width=9cm]{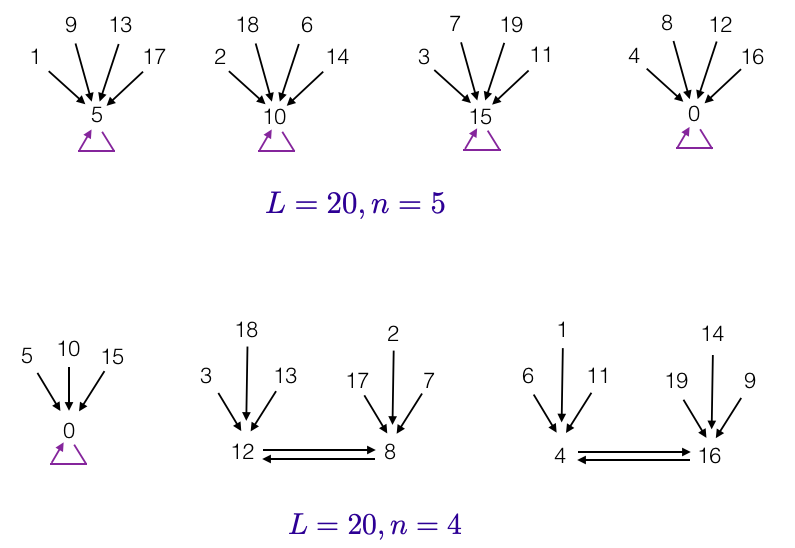}
\caption{ Bottom: $Z(4, 20)$ has $15$ leaves, $5$ roots, $0$ branches, and $3$ connected components.  Top: $Z(5, 20)$ has $16$ leaves, $4$ roots, $0$ branches, and  $4$ connected components. } 
\label{newL=20,n=4,5}
\end{figure}

\begin{figure}[H]
\includegraphics[width=9cm]{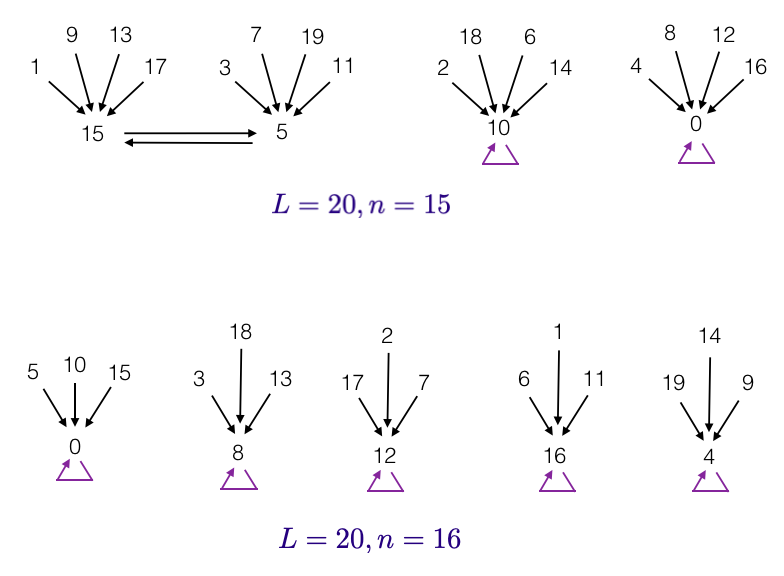}
\caption{Top: $Z(15, 20)$ has $16$ leaves, $4$ roots, $0$ branches, and $3$ connected components.
Bottom: $Z(16, 20)$ has $15$ leaves, $5$ roots, $0$ branches, and $5$ connected components.} 
\label{new,L=20,n=15,16}
\end{figure}

\begin{figure}[H]
\includegraphics[width=8cm]{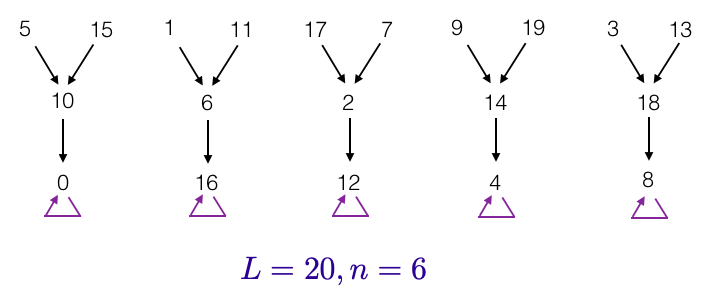}
\caption{$Z(6,20)$ has $10$ leaves, $5$ branches, $5$ roots, and $5$ connected components.} 
\label{newL=20,n=6}
\end{figure}

\begin{figure}[ht]
\includegraphics[width=4cm]{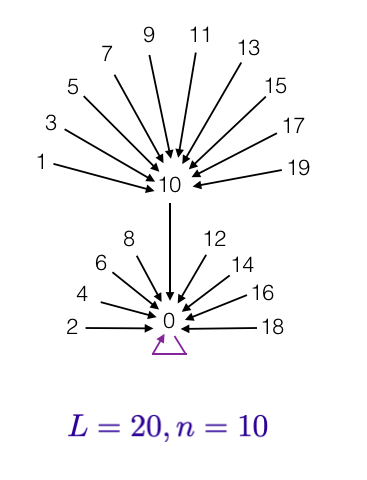}
\caption{$Z(10, 20)$ has $18$ leaves, $1$ root, $1$ branch, and just $1$ connected component.} 
\label{second appearance of Z(10,20)}
\end{figure}

 The reader may consult \url{https://github.com/Caio-Simon/Zolotarev-Graphs} for code that constructs Zolotarev graphs for any $L$ and $n$.

\bigskip
%%%%%%%%%%%%%%%%%%%%%%%%%%%
%%%%%%%%%%%%%%%%%%%%%%%%%%%%%%%

\end{document}